\documentclass[%
 aip,
 amsmath,amssymb,
reprint,%
]{revtex4-1}

\usepackage{graphicx}
\usepackage{dcolumn}
\usepackage{bm}

\usepackage[utf8]{inputenc}
\usepackage[T1]{fontenc}
\usepackage{mathptmx}
\usepackage{etoolbox}

\usepackage{hyperref}
\hypersetup{
    colorlinks=true,
    linkcolor=blue,
    filecolor=magenta,      
    urlcolor=blue,
    citecolor=blue,
}

\makeatletter
\def\@email#1#2{%
 \endgroup
 \patchcmd{\titleblock@produce}
  {\frontmatter@RRAPformat}
  {\frontmatter@RRAPformat{\produce@RRAP{*#1\href{mailto:#2}{#2}}}\frontmatter@RRAPformat}
  {}{}
}%
\makeatother
\begin{document}

\preprint{AIP/123-QED}

\title[Invariant Curves and the Variational Structure in Tubular Origami Dynamical Systems]{Invariant Curves and the Variational Structure in Tubular Origami Dynamical Systems}

\author{Ryutaro. Ichikawa}
\author{Mitsuru. Shibayama}%
 \email{ichikawa.ryutaro.y81@kyoto-u.jp and shibayama@amp.i.kyoto-u.ac.jp}
\affiliation{ 
Dynamical Systems Group,
Applied Mathematics and Physics Course,
Graduate School of Informatics,
Kyoto University,
Yoshida-honmachi, Sakyo-ku, Kyoto, 606-8501, JAPAN
}%

\date{\today}

\begin{abstract}
    We present a theoretical and numerical dynamical-systems analysis of tubular origami tessellations by identifying the inverse module number, $N^{-1}$, as a perturbation parameter within the framework of Kolmogorov--Arnold--Moser (KAM) theory. In the large-module limit ($N \to \infty$), we show that the conservative dynamics formally converges to an integrable map with a variational structure, whose generating function corresponds to the total discrete mean curvature.
    From the viewpoint of KAM theory, nonresonant invariant curves of the integrable limit are expected to persist for sufficiently large \(N\).
    Consistent with this expectation, numerical computations with increasing \(N\) show that large regions of the phase space are filled with structures that appear to be invariant curves. By adjusting mountain-valley fold assignments and fold lengths, the system can be transformed into a nontwist map that exhibits multiple zero frequencies. 
    The frequency profile in the integrable limit and the persistence of invariant curves allow us to control the number and arrangement of stable folding regions appearing as coexisting elliptic islands. 
    These islands provide a phase-space interpretation of distinct folding modes, separated by invariant curves that act as geometric barriers to continuous deformation and obstruct transitions without self-intersection.
    Finally, we analyze the expanding and contracting dynamics of the origami structure within the framework of conformally symplectic systems. By introducing a virtual auxiliary fold as a drift control mechanism, we numerically confirm the existence of stable quasi-periodic attractors.
\end{abstract}

\maketitle

\begin{quotation}
    \emph{Origami tessellations}, consisting of periodic crease patterns, exhibit highly nontrivial nonlinear folding behavior when individual modules deform non-uniformly. In the case of tubular origami tessellations, the spatial propagation of fold states can be formulated as an \emph{area-preserving discrete dynamical map}. Although these maps admit a variational formulation through generating functions, the geometric meaning of these generating functions is difficult to identify explicitly because of the analytical complexity of the folding map.
    In this study, we analyze a formal large-module limit in which the number of modules tends to infinity. This integrable limit clarifies the \emph{variational structure} of the folding map. Specifically, we show that the associated \emph{generating function} is identical to a fundamental geometric quantity, namely, \emph{the total discrete mean curvature} of the folded surface. This result reveals a direct connection between the discrete geometry of folded surfaces and the variational structure underlying tubular origami dynamics.
    The integrable-limit analysis is also useful for interpreting the phase-space organization of the finite-\(N\) system. 
    In particular, invariant curves associated with self-intersecting folded shapes can be interpreted as geometric barriers to continuous deformation between stable folding regions.
\end{quotation}
    
\section{\label{sec:intro}INTRODUCTION}
    \begin{figure*}
        \includegraphics[width=15cm]{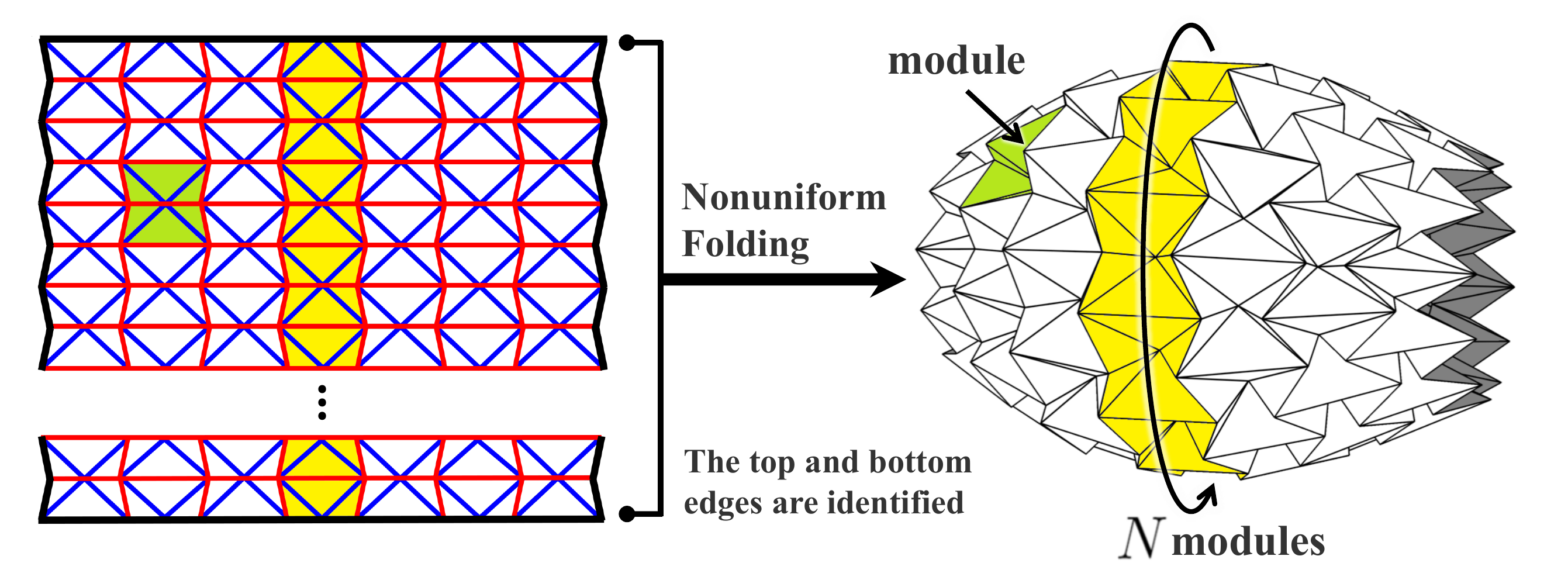}
        \caption{\label{fig:origami_model}
            Schematic illustration of the waterbomb tube origami tessellation.
            Left: The crease pattern where red, blue, and black lines represent mountain folds, valley folds, and boundary lines, respectively. The top and bottom edges are identified. Right: The 3D cylindrical structure obtained by nonuniform-folding. The light green polygon indicates a single module, and the yellow region highlights a ring consisting of $N$ modules. These colored regions correspond to each other in both figures.
        }
    \end{figure*}

    Origami-inspired mechanical metamaterials have garnered significant interest as adaptable structural systems because of their high deployability and adjustable mechanical stiffness. Notable applications span a wide range of fields\cite{comprehensive_review}, including deployable solar panels for aerospace systems\cite{Yue_2023}, biomedical stents\cite{KURIBAYASHI2006131}, architectural structures\cite{Vlachaki2025}, and soft robotics\cite{SoftRobots}.
    
    Origami structures are also known to exhibit rich nonlinear phenomena associated with multistability.
    Multiple energy-stable configurations can coexist in a single origami structure, and transitions between these states can be induced by external loading, snap-through instabilities, or dynamic excitation\cite{PhysRevE.95.052211,SadeghiLi2021}.
    Such multistable responses have been exploited to design tunable origami units, mechanical memory devices, and programmable metamaterials\cite{YasudaTachiLeeYang2017,LiuTachiPaulino2019}.
    Recent studies have further investigated complex transition sequences and nonlinear dynamic responses in multistable origami structures, including intra-well and inter-well motions and data-driven identification of governing equations\cite{PhysRevE.95.052211,SadeghiLi2021,LIU2024111659}.
    These nonlinear responses originate from the geometric constraints of origami.
    However, when such constraints are considered in tessellated structures with many coupled modules, direct analysis becomes difficult.
    
    Owing to this difficulty, 
    much of existing engineering-oriented research on origami tessellations has been built on uniformly folded configurations, where all unit cells deform identically. While such \emph{uniform-folding} provides a convenient foundation for understanding and controlling the kinematics and mechanical responses of periodic crease patterns\cite{uniform0, uniform1, uniform2, uniform3}, 
    the range of achievable mechanical responses and geometric configurations remains inherently limited. 
    To overcome these limitations, 
    it is necessary to understand \emph{nonuniform-folding} processes, in which individual modules deform differently.
    Several studies have investigated nonuniform-folding processes in origami tessellations from geometric and mechanical perspectives\cite{nonuniform0,nonuniform1,ZHAO2018442}.
    In such processes, the folding state of one module constrains the state of neighboring modules through the geometry of the crease pattern.
    This makes the spatial propagation of folding states a central object of study.
    Complementary to energy-based and excitation-driven studies of origami nonlinearity, recent geometric approaches have formulated this propagation process as a discrete dynamical system derived from the kinematic constraints of origami\cite{ImadaTachi2022,ImadaTachi2023,ImadaTachi2025}.
    This viewpoint describes folding evolution not as temporal motion under external forcing, but as spatial evolution along the tessellation.
    In particular, for tubular tessellations, it has been established that this spatial propagation is governed by a \emph{symplectic (area-preserving)} map\cite{ImadaTachi2023}, indicating the existence of conserved quantities or invariant structures analogous to those found in Hamiltonian mechanics\cite{RevModPhys.64.795}.

    However, the resulting symplectic map has a highly complex analytical structure, making it difficult to explicitly identify the underlying variational structure. To address this challenge, we consider the regime with a sufficiently large number of modules. In this limit, the system approaches an integrable state, where the variational structure becomes analytically accessible. This integrable system reveals the geometric origin of the conservation laws and provides a baseline for understanding the global dynamics, as the persistence of invariant curves is supported theoretically by \emph{KAM theory}\cite{Lichtenberg1992, Arnold1989}.
    
    In this study, we conduct a global dynamical systems analysis of origami tessellations and examine their variational structure by treating the inverse module number, $N^{-1}$, as a \emph{perturbation} parameter within the framework of KAM theory. We demonstrate that, when the number of modules is sufficiently large, the system simplifies to an integrable dynamical system and derive its generating function. Specifically, for symmetric crease patterns, we identify this generating function as the total discrete mean curvature, providing a clear geometric interpretation of the conserved quantity that governs the dynamics of tubular origami.
    
    Based on this integrable structure, we establish the existence of invariant curves and confirm their persistence under finite-$N$ perturbations, as predicted by KAM theory. Furthermore, our analysis of the generating function reveals that the system's frequency depends sensitively on the \emph{mountain–valley} fold assignments. By appropriately tuning these parameters, we can violate the \emph{twist condition}, leading to the emergence of \emph{nontwist} dynamical behavior in tubular origami tessellations. In the nontwist regime, the frequency profile can exhibit multiple zero values, and the resonances associated with these points give rise to new stable foldable regions in phase space, appearing as islands of elliptic periodic orbits. We numerically confirm the emergence of these structures and establish their correspondence with physically realizable folding configurations. Previous studies\cite{ImadaTachi2023} have interpreted quasiperiodic orbits around elliptic fixed points as nonuniform folded configurations, including undulated tubular shapes. 
    Within this phase-space description, a continuous displacement between nearby orbits corresponds to a continuous deformation of the associated folded shapes. 
    Our numerical computations reveal that some invariant curves separating stable islands are associated with self-intersecting folded shapes. 
    This observation provides a geometric interpretation of these invariant curves: they act as barriers to continuous deformation between distinct stable folding regions.
    
    Finally, we observe that the two-degree-of-freedom origami system can be formulated as a \emph{conformally symplectic (CS) system}, characterized by a constant Jacobian determinant. It is known that CS systems can support stable \emph{quasi-periodic attractors} through the introduction of a drift term that counteracts phase-space expansion or contraction\cite{CALLEJA2013978}. Motivated by this insight, we introduce a \emph{virtual fold} designed to generate such a \emph{drift term} and numerically demonstrate the existence of stable quasi-periodic attractors in the resulting folding dynamics.
    
    This paper is organized as follows. Sec.~\ref{sec:model} presents the tubular origami model and describes the discrete dynamical system that governs its folding behavior, emphasizing its CS structure. Sec.~\ref{sec:theory} describes the theoretical analysis, demonstrating the reduction to an integrable system, clarifying the geometric interpretation of the underlying variational structure, and examining the existence of invariant curves using generating functions. Sec.~\ref{sec:numerical} presents the results of numerical experiments, highlighting phase space structures such as KAM tori and a KAM attractor. Finally, Sec.~\ref{sec:conclusion} summarizes our main conclusions.

\section{\label{sec:model}THE TUBULAR ORIGAMI MODEL}

    In this section, we introduce the discrete dynamical system for tubular origami tessellations proposed in previous studies\cite{ImadaTachi2023}, which serves as the fundamental framework throughout this paper. As illustrated in Fig.~\ref{fig:origami_model}, we consider an origami tessellation composed of identical modules arranged in a periodic crease pattern (left), which forms a three-dimensional cylindrical structure upon folding (right). A collection of $N$ modules ($N \in \mathbb{Z}_{>2}$), corresponding to the yellow region in the crease pattern, is referred to as a \emph{ring}. This ring serves as the basic unit of the discrete dynamical system considered in this study.
    
    \begin{figure}
        \includegraphics[width=\columnwidth]{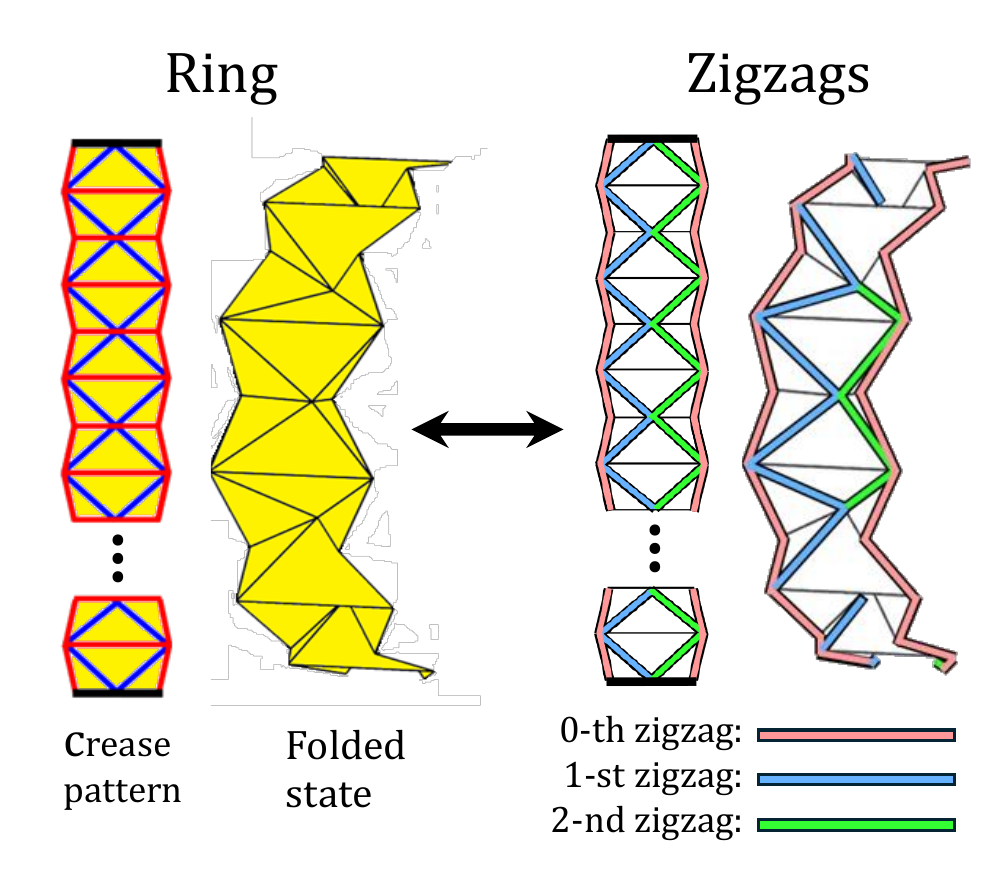}
        \caption{\label{fig:zigzag_model}
            Ring and zigzags as constitutive elements. (Left) The crease pattern of the ring shown in Fig.~\ref{fig:origami_model} and its corresponding folded configuration. (Right) Structural decomposition into zigzags. The solid lines in pink, light blue, and yellow-green represent adjacent zigzags (the 0-th, 1-st, and 2-nd zigzags, respectively), illustrating how the global structure is formed by their concatenation along the cylindrical axis.
        }
    \end{figure}

\subsection{\label{sec:param}Modeling and Parameterization}
    As a representative example, we consider the waterbomb tube, which serves as a canonical model for tubular origami tessellations. As illustrated in Fig.~\ref{fig:zigzag_model}, the entire structure can be decomposed into a periodic sequence of closed polygonal chains, referred to as \emph{zigzags}. We focus on the discrete evolution of the zigzag geometry along the cylindrical axis. By assuming an $N$-fold rotational symmetry of the folded configuration, the geometric information of an entire ring is reduced to that of a single module. Accordingly, we introduce a parametrization of the zigzag associated with a single unit cell.
        
\begin{figure}
        \includegraphics[width=\columnwidth]{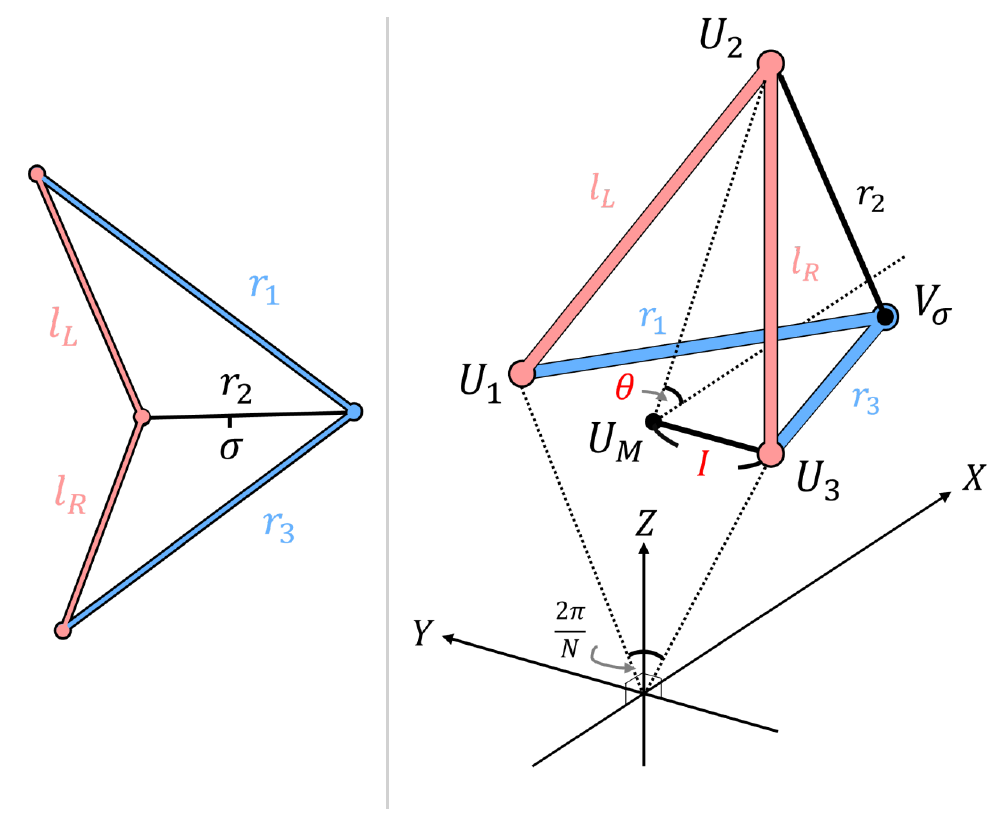}
        \caption{\label{fig:parametrization_zigzag}
            Parametrization of the zigzag within a module. The zigzag segment defined by vertices $U_1, U_2$, and $U_3$ (the $i$-th zigzag) is uniquely determined by the state variables $(\theta, I)$. Here, $I$ is defined as half the Euclidean distance between $U_1$ and $U_3$, and $\theta$ is defined as the angle formed by the vector $\overrightarrow{U_M U_2}$ relative to the central cylinder axis ($X$-axis). The position of the vertex $V_\sigma$ in the subsequent ($i+1$)-th zigzag is determined by the intersection of three spheres centered at $U_1, U_2$, and $U_3$.
        }
    \end{figure}

    First, we define a global Cartesian coordinate system $(X,Y,Z)$ such that the $X$-axis aligns with the central axis of the folded cylinder. Within this framework, the continuous zigzag path wrapping around the cylinder can be interpreted as a concatenation of discrete zigzag segments derived from individual modules, each obtained by successive rotation about the $X$-axis by an angle $2\pi/N$. Focusing on a single representative module as shown in Fig.~\ref{fig:parametrization_zigzag}, we denote the vertices of the zigzag element by $U_1$, $U_2$, and $U_3$. To establish a local reference frame for this module, the $Y$-axis is chosen to align with the direction of the vector $\overrightarrow{U_3U_1}$. The $Z$-axis is defined to be vertical such that the endpoints $U_1$ and $U_3$ lie in the $YZ$-plane (i.e., the plane $X=0$).
        
    Let $\mathbf{r}_{U_j}$ $(j=1,2,3)$ denote the position vectors of the vertices $U_j$. The geometry of the zigzag element is characterized by prescribed positive crease lengths $l_L \in \mathbb{R}_{>0}$ and $l_R \in \mathbb{R}_{>0}$, which impose geometric distance constraints on the vertex positions. We define $\mathbf{r}_{V_\sigma}$ as a point determined by the intersection of three spheres centered at $\mathbf{r}_{U_j}$, each with radius $r_j>0$ $(j=1,2,3)$, namely,
    \[
        \|\mathbf{r}_{V_\sigma}-\mathbf{r}_{U_j}\| = r_j, \qquad j=1,2,3.
    \]
    Since the intersection of three spheres generically yields two distinct solutions, the physically realized configuration is selected according to the mountain--valley assignment $\sigma$.
    
    Using this geometric setup, we introduce a set of state variables $(\theta,I) \in (-\pi,\pi] \times \left(|l_L-l_R|/2, (l_L+l_R)/2\right)$ to characterize the geometry of a single zigzag module. These variables are defined by
    \begin{eqnarray}
        I &=& \frac{1}{2}\left\| \mathbf{r}_{U_3} - \mathbf{r}_{U_1} \right\|, \label{eq:def_I} \\
        \theta &=& \arctan2\!\left(
        \left(\mathbf{e}_X\times \frac{\mathbf{v}}{\|\mathbf{v}\|}\right)\cdot 
        \frac{\mathbf{r}_{U_3}-\mathbf{r}_{U_1}}{2I},
        \mathbf{e}_X \cdot \frac{\mathbf{v}}{\|\mathbf{v}\|}
        \right), \label{eq:def_theta}
    \end{eqnarray}
    where $\mathbf{e}_X = [1, 0, 0]^T$ and $\mathbf{v} = \mathbf{r}_{U_2} - \mathbf{r}_{U_M}$. Here, $\mathbf{r}_{U_M}$ denotes the foot of the perpendicular dropped from $U_2$ onto the line segment $U_1U_3$.
        
    The explicit expressions for the position vectors $\mathbf{r}_{U_1}$, $\mathbf{r}_{U_2}$, $\mathbf{r}_{U_3}$, and $\mathbf{r}_{U_M}$ are given by
    \begin{eqnarray}
        \mathbf{r}_{U_1} &=& \left[0, I, I\cot\frac{\pi}{N} \right]^T, \\
        \mathbf{r}_{U_2} &=& \left[0, I-l_L\alpha(I), I\cot\frac{\pi}{N} \right]^T \nonumber\\ 
        & \ &+ l_L\beta(I)\left[\cos\theta, 0, \sin\theta \right]^T, \\
        \mathbf{r}_{U_3} &=& \mathbf{R}_X\!\left(\frac{2\pi}{N}\right)\mathbf{r}_{U_1}
        = \left[0, -I, I\cot\frac{\pi}{N} \right]^T, \\
        \mathbf{r}_{U_M} &=& \mathbf{r}_{U_1} 
        + l_L\alpha(I)\left(\frac{\mathbf{r}_{U_3}-\mathbf{r}_{U_1}}{2I}\right),
    \end{eqnarray}
    where $\mathbf{R}_X\!\left(2\pi/N\right)$ denotes the rotation matrix associated with a rotation of angle $2\pi/N$ about the $X$-axis. The functions $\alpha(I)$ and $\beta(I)$ are defined as
    \begin{equation}
        \alpha(I) = \frac{l_L^2 + (2I)^2 - l_R^2}{2l_L(2I)},
        \qquad
        \beta(I) = \sqrt{1 - \alpha^2(I)}.
    \end{equation}
    
    Next, we introduce a local coordinate system attached to the vertex $U_1$. Let the origin be located at $\mathbf{r}_{U_1}$, and define an orthonormal basis $\{\mathbf{e}_1,\mathbf{e}_2,\mathbf{e}_3\}$ by
    \begin{eqnarray}
        \mathbf{e}_1 &=& \frac{\mathbf{r}_{U_2}-\mathbf{r}_{U_1}}
        {\|\mathbf{r}_{U_2}-\mathbf{r}_{U_1}\|} = \left[\beta\cos\theta, -\alpha, \beta\sin\theta\right]^T, \label{eq:e1_def} \\
        \mathbf{e}_2 &=&
        \frac{(\mathbf{r}_{U_3}-\mathbf{r}_{U_1})
        - \left((\mathbf{r}_{U_3}-\mathbf{r}_{U_1})\cdot\mathbf{e}_1\right)\mathbf{e}_1}
        {\sqrt{\|\mathbf{r}_{U_3}-\mathbf{r}_{U_1}\|^2
        - \left((\mathbf{r}_{U_3}-\mathbf{r}_{U_1})\cdot\mathbf{e}_1\right)^2}} \label{eq:e2_def} \nonumber \\
        &=& -\left[\alpha\cos\theta, \beta, \alpha\sin\theta\right]^T \\
        \mathbf{e}_3 &=& \mathbf{e}_1 \times \mathbf{e}_2 = \left[\sin\theta, 0, -\cos\theta\right]^T. \label{eq:e3_def}
    \end{eqnarray}
        
    In this local coordinate system, the position vector of $V_\sigma$ is expressed as
    \begin{equation}
        \mathbf{r}_{V_\sigma}
        = \mathbf{r}_{U_1} + e_1 \mathbf{e}_1 + e_2 \mathbf{e}_2 + e_3 \mathbf{e}_3,
    \end{equation}
    where the scalar components $e_1$, $e_2$, and $e_3$ are determined by the intersection of three spheres centered at $\mathbf{r}_{U_1}$, $\mathbf{r}_{U_2}$, and $\mathbf{r}_{U_3}$, with radii $r_1>0$, $r_2>0$, and $r_3>0$, respectively. The components $e_1$, $e_2$, and $e_3$ are explicitly given by
    \begin{eqnarray}
        e_1 &=&
        \frac{r_1^2 + \|\mathbf{r}_{U_2}-\mathbf{r}_{U_1}\|^2 - r_2^2}
        {2\|\mathbf{r}_{U_2}-\mathbf{r}_{U_1}\|}=\frac{r_1^2+l_L^2-r_2^2}{2l_L}, \label{eq:e1} \\
        e_2 &=&
        \frac{1}{2\,(\mathbf{r}_{U_3}-\mathbf{r}_{U_1})\cdot\mathbf{e}_2}
        \Bigl(
        r_1^2 - r_3^2
        + \left((\mathbf{r}_{U_3}-\mathbf{r}_{U_1})\cdot\mathbf{e}_1\right)^2
        \nonumber\\
        &&
        + \left((\mathbf{r}_{U_3}-\mathbf{r}_{U_1})\cdot\mathbf{e}_2\right)^2
        - 2e_1(\mathbf{r}_{U_3}-\mathbf{r}_{U_1})\cdot\mathbf{e}_1
        \Bigr) \nonumber \\
        &=& \frac{r_1^2-r_3^2+4I^2 -4Ie_1\alpha}{4I\beta}
        \label{eq:e2} \\
        e_3 &=& \sigma\sqrt{\gamma}, \label{eq:e3}
    \end{eqnarray}
    where $\gamma = r_1^2 - e_1^2 - e_2^2$ and the sign $\sigma=\pm1$ corresponds to the mountain--valley assignment ($\sigma(M)=1, \sigma(V)=-1$). 
        
    Finally, as illustrated in Fig.~\ref{fig:zigzag_model}, the $(i+1)$-th zigzag is generated from the $i$-th zigzag with a one-step circumferential offset. Accordingly, when computing the parameters $(\theta', I')$ of the next zigzag, we select the vertices of the zigzag element as follows: the position vector of the new $U_1$ is taken to be $\mathbf{r}_{V_\sigma}$, that of the new $U_2$ is chosen as $\mathbf{r}_{U_3}$, and that of the new $U_3$ is given by $\mathbf{R}_X\!\left(2\pi/N\right)\mathbf{r}_{V_\sigma}$. Applying the definitions of $I$ and $\theta$ given in Eqs.~\eqref{eq:def_I} and \eqref{eq:def_theta} to this new triplet of vertex position vectors yields the updated state variables $(\theta', I')$. This construction defines a discrete evolution of the zigzag geometry along the axial direction of the tubular origami.
        
    \subsection{\label{sec:dynsys}Definition of the Discrete Dynamical System}
    
        Based on the geometric construction and parametrization introduced in Sec.~\ref{sec:param}, we now define the discrete dynamical system that governs the evolution of the zigzag geometry along the axial direction of the tubular origami. 
        Let the map $\mathbf{f}$ that sends $(\theta, I)$ to $(\theta^{\prime}, I^{\prime})$ through the three-sphere intersection be defined as
        \begin{eqnarray}
            \lefteqn{\left[\theta^{\prime}, I^{\prime}\right]^T := \mathbf{f}\left(\left[\theta, I\right]^T;l_L, l_R, N, r_1, r_2, r_3, \sigma\right)} \label{eq:map_f_def} \\
            \theta^{\prime} &=& \arctan2\left(
                \left(\mathbf{e}_X\times\frac{\mathbf{v}}{\|\mathbf{v}\|}\right)\cdot \frac{\mathbf{u}}{\|\mathbf{u}\|}, 
                \mathbf{e}_X \cdot \frac{\mathbf{v}}{\|\mathbf{v}\|}
            \right), \label{eq:map_f_theta} \\ 
            I^{\prime} &=& \frac{1}{2} \|\mathbf{u}\| \label{eq:map_f_I},
        \end{eqnarray}
        where $\mathbf{E}$ denotes the $3 \times 3$ identity matrix and
        \begin{eqnarray}
            \mathbf{u} &=& 
            \begin{bmatrix}
            0 \\ -2I \\ 0
            \end{bmatrix}
            + \left(\mathbf{R}_X\!\left(\frac{2\pi}{N}\right) - \mathbf{E}\right)
            \left(\mathbf{e}_1 \ \mathbf{e}_2 \ \mathbf{e}_3 \right)
            \begin{bmatrix}
            e_1 \\ e_2 \\ e_3
            \end{bmatrix}, \\
            \mathbf{v} &=& \mathbf{r}_{U_3} - \left(
                \mathbf{r}_{V_\sigma} + r_3\left(
                \frac{r_3^2+(2I')^2-r_1^2}{2r_3(2I')}
                \right)
                \frac{\mathbf{u}}{2I'}
            \right).
        \end{eqnarray}
        When the tube consists of $m \in \mathbb{Z}_{>0}$ zigzags indexed by $i=0, \dots, {m-1}$, we allow the geometric parameters to vary along the axial direction.
        For the $i$-th zigzag, the left, middle, and right crease lengths are denoted by $l_{i,L}$, $l_{i,M}$, and $l_{i,R}$, respectively.
        To ensure consistency with the periodic arrangement of identical modules along the axial direction, we impose periodic boundary conditions such that $l_{m,L} = l_{0,L}$, $l_{m,M} = l_{0,M}$, and $l_{m,R} = l_{0,R}$.
        The map $\mathbf{f}_i$ describing the transition from the $i$-th zigzag to the $(i+1)$-th zigzag is then defined as follows:
        \begin{equation}
            \mathbf{f}_i\left(\left[\theta, I\right]^T\right) = \mathbf{f}\left(\left[\theta, I\right]^T;l_{i,L}, l_{i,R}, N, l_{i+1, R}, l_{i+1,M}, l_{i+1, L}, \sigma_i\right).
        \end{equation}

        \begin{figure}
            \includegraphics[width=\columnwidth]{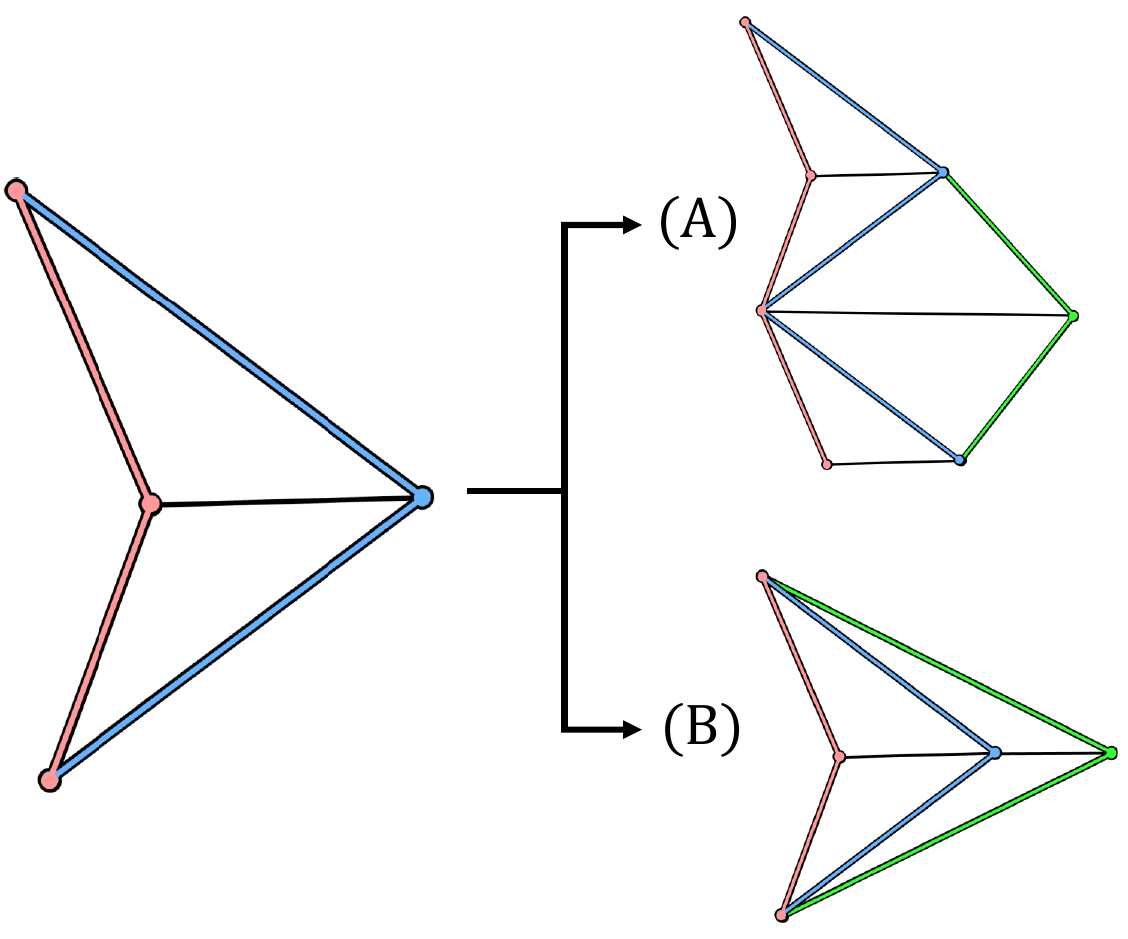}
            \caption{\label{fig:2type_zigzag}
                Two distinct patterns of zigzag connectivity. (A) The configuration involving a circumferential shift, which corresponds to the waterbomb tube structure. (B) The configuration where $V_\sigma$ is constructed directly from the original zigzag without a shift. In this pattern, a re-parametrization of the state variables $(\theta, I)$ is required.
            }
        \end{figure}
        
        \begin{figure*}
            \includegraphics[width=18cm]{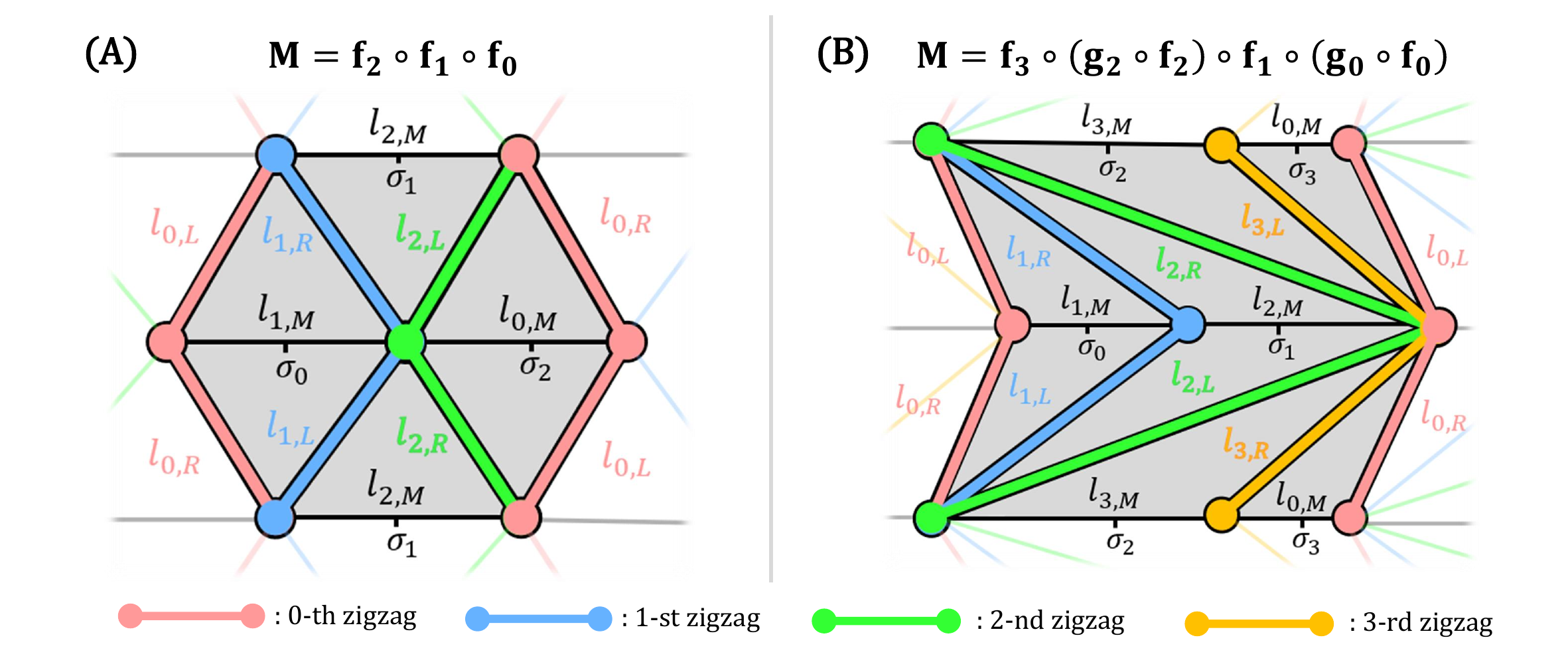}
            \caption{\label{fig:examples_of_modules}
                Crease patterns for two representative module configurations. (A) A configuration featuring six-valent vertices, which corresponds to a generalization of the waterbomb tessellation. (B) A configuration characterized by four- or eight-valent vertices. In both panels, different colors indicate distinct zigzags, and the gray shaded region represents the fundamental domain of a single module.
            }
        \end{figure*}
        
        As shown in Fig.~\ref{fig:2type_zigzag}, zigzag modules can be connected in two topologically distinct ways. The map $\mathbf{f}$ considers only type~(A) connectivity, where consecutive zigzags are joined with a circumferential offset. For modules with type~(B) connectivity, however, the parameters $(\theta, I)$ obtained from $\mathbf{f}$ do not directly represent the geometric state of the next zigzag and must be reparametrized. To account for this, we introduce a  map $\mathbf{g}$:
        \begin{eqnarray}
            \left[\theta^{*}, I^{*}\right]^T &:=& \mathbf{g}\left(\left[\theta, I\right]^T;r_3, r_1, N\right), \label{eq:map_g_def} \\
            \theta^{*} &=& \arctan2\left(
                Y_g, X_g
            \right), \label{eq:map_g_theta} \\ 
            I^{*} &=& \frac{\sin\frac{\pi}{N}}{4I}\sqrt{\left(r_3^2-r_1^2\right)^2 + T_2^2} \label{eq:map_g_I} ,
        \end{eqnarray}
        where 
        \begin{eqnarray}
            T_1 &=& \sqrt{16I^2 r_3^2 - \left(4I^2 + r_3^2 - r_1^2\right)^2}, \\ 
            T_2 &=& 4I^2\cot\frac{\pi}{N} + T_1\sin\theta, \\
            Y_g &=& -\sin\frac{\pi}{N}\Biggl(
            T_2\left(T_1\cot\frac{\pi}{N}\sin\theta -4I^2\right) \nonumber \\
              && \qquad+ \left(r_1^2-r_3^2\right)^2\cot\frac{\pi}{N}
            \Biggr), \label{eq:def_Yg} \\
            X_g &=& -T_1\cos\theta \sqrt{\left(r_3^2-r_1^2\right)^2 + T_2^2} . \label{eq:def_Xg}
        \end{eqnarray}
        
        Here, we introduce a binary symbol $k_i \in \{0,1\}$ to indicate whether the re-parametrization map $\mathbf{g}$ is applied at the $i$-th zigzag.
        Accordingly, we define the switching map $\mathbf{g}_i^{k_i}$ as
        \begin{equation}
            \mathbf{g}_i^{k_i}\left(\left[\theta, I\right]^T\right) :=
            \begin{cases}
                \left[\theta, I\right]^T & \mathrm{if} \ \ k_i = 0, \\
                \mathbf{g}\left(\left[\theta, I\right]^T ; l_{i+1, L}, l_{i+1, R}, N \right) & \mathrm{if} \ \ k_i = 1.
            \end{cases}
        \end{equation}
        In applying the map $\mathbf{g}$, certain operations such as the interchange of $l_L$ and $l_R$ are required. The details of these operations are not essential for the present discussion and are therefore deferred to previous work\cite{ImadaTachi2023}.
        Finally, for each zigzag indexed by $i$, we define the corresponding state transition map $\mathbf{M}_i$ as the composition of $\mathbf{f}_i$ and $\mathbf{g}_i^{k_i}$:
        \begin{equation}
            \mathbf{M}_i
            :=
            \mathbf{g}_i^{\,k_i}
            \circ
            \mathbf{f}_i  \qquad (i=0, \ldots, m-1).
        \end{equation}
        The evolution of the state variables $(\theta, I)$ across a ring composed of $m$ zigzags is determined by the composition of a sequence of maps $\mathbf{M}_i$.
       By defining the global transfer map $\mathbf{M}$ as this composition, the discrete dynamical system that governs the geometric evolution from the $n$-th ring to the $(n+1)$-th ring is formulated as:
         \begin{eqnarray}
            \mathbf{M} &:=& \mathbf{M}_{m-1} \circ \cdots \circ \mathbf{M}_0, \label{eq:map_M_def} \\
            \left[\theta_{n+1}, I_{n+1}\right]^T &=& \mathbf{M}\left( \left[\theta_{n}, I_{n}\right]^T \right). \label{eq:dynamical_system}
        \end{eqnarray}
        It is important to note that this discrete dynamical system is area-preserving, meaning that the map $\mathbf{M}$ preserves the symplectic structure.

    \subsection{\label{sec:example_of_modules}Examples of Conservative Discrete Dynamical Systems}
        
        Following the formulation of the conservative discrete dynamical system in the previous section, we now present two representative examples of module configurations. Shown in Fig.~\ref{fig:examples_of_modules}, these examples display distinct connectivity patterns and valencies.
        
        In Case (A) (see Fig.~\ref{fig:examples_of_modules}), the module consists of three zigzags ($m=3$) connected only through type-(A) connectivity (see Fig.~\ref{fig:2type_zigzag}), which implies $k_i=0$ for all $i$.
        Consequently, the return map $\mathbf{M}$ is defined as the composition $\mathbf{f}_2 \circ \mathbf{f}_1 \circ \mathbf{f}_0$, where each map is given explicitly by
        \begin{eqnarray}
            \mathbf{f}_0(\mathbf{x}) &=& \mathbf{f}\left(\mathbf{x}; l_{0,L}, l_{0,R}, N, l_{1,R}, l_{1,M}, l_{1,L}, \sigma_0\right), \\
            \mathbf{f}_1(\mathbf{x}) &=& \mathbf{f}\left(\mathbf{x}; l_{1,L}, l_{1,R}, N, l_{2,R}, l_{2,M}, l_{2,L}, \sigma_1\right), \\
            \mathbf{f}_2(\mathbf{x}) &=& \mathbf{f}\left(\mathbf{x}; l_{2,L}, l_{2,R}, N, l_{0,R}, l_{0,M}, l_{0L}, \sigma_2\right),
        \end{eqnarray}
        where $\mathbf{x} = [\theta, I]^T$.
        In this configuration, the absence of re-parametrization ($\mathbf{g}_i$) preserves the six-valent vertex structure throughout the module, maintaining the standard topological connectivity.
        
        In contrast, Case (B) (see Fig.~\ref{fig:examples_of_modules}) represents a module with four zigzags ($m=4$) incorporating type-(B) connectivity (refer to Fig.~\ref{fig:2type_zigzag}).
        As indicated in Fig.~\ref{fig:examples_of_modules}(B), the re-parametrization map $\mathbf{g}$ is applied at the 0-th and 2-nd zigzags ($k_0=1, k_2=1$), while it is omitted at the 1-st and 3-rd zigzags ($k_1=0, k_3=0$).
        Consequently, the return map is expressed as $\mathbf{M} = \mathbf{f}_3 \circ (\mathbf{g}_2 \circ \mathbf{f}_2) \circ \mathbf{f}_1 \circ (\mathbf{g}_0 \circ \mathbf{f}_0)$, where the constituent maps are defined as follows:
        \begin{eqnarray}
            \mathbf{f}_0(\mathbf{x}) &=& \mathbf{f}\left(\mathbf{x}; l_{0,L}, l_{0,R}, N, l_{1,R}, l_{1,M}, l_{1,L}, \sigma_0\right), \\
            \mathbf{g}_0(\mathbf{x}) &=& \mathbf{g}\left(\mathbf{x}; l_{1,L}, l_{1,R}, N\right), \\
            \mathbf{f}_1(\mathbf{x}) &=& \mathbf{f}\left(\mathbf{x}; l_{1,R}, l_{1,L}, N, l_{2,R}, l_{2,M}, l_{2,L}, \sigma_1\right), \\
            \mathbf{f}_2(\mathbf{x}) &=& \mathbf{f}\left(\mathbf{x}; l_{2,L}, l_{2,R}, N, l_{3,R}, l_{3,M}, l_{3,L}, \sigma_2\right), \\
            \mathbf{g}_2(\mathbf{x}) &=& \mathbf{g}\left(\mathbf{x}; l_{3,L}, l_{3,R}, N\right), \\
            \mathbf{f}_3(\mathbf{x}) &=& \mathbf{f}\left(\mathbf{x}; l_{3,R}, l_{3,L}, N, l_{0,R}, l_{0,M}, l_{0,L}, \sigma_3\right).
        \end{eqnarray}
        This repeated application of the re-parametrization map modifies the topological structure, producing a tessellation with vertices of valence four or eight.
    
    \subsection{\label{sec:cs}CS Structure and Drift}
        In addition to the conservative systems discussed in the preceding section, the previous study\cite{ImadaTachi2023} also examines tessellations that allow for similarity expansion. Specifically, this is achieved by introducing a scaling factor $s \in \mathbb{R}_{>0}$, so that if the crease lengths of the initial zigzag are $l_{0,L}$ and $l_{0,R}$, the lengths of the final zigzag become $s l_{0,L}$ and $s l_{0,R}$. However, applying the scaling factor $s$ also scales the domain of the state variables. To construct a consistent dynamical system in which the phase space remains invariant across iterations, the output domain must be aligned with the initial domain. This is accomplished by multiplying the output by $1/s$ at the final stage of the mapping. As a result, the determinant of the Jacobian matrix of the system becomes $1/s$, making it a dissipative system rather than a conservative one.
        
\begin{figure}
            \includegraphics[width=\columnwidth]{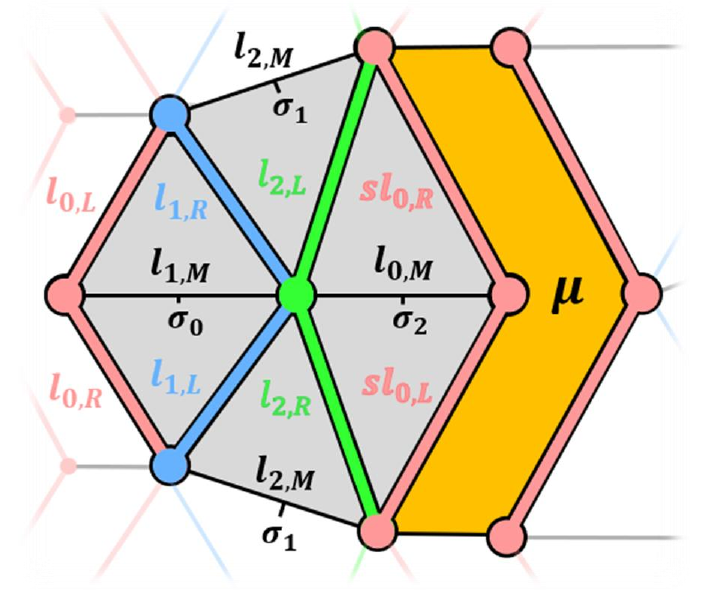}
            \caption{\label{fig:cs_pattern} Crease pattern of the extended module designed for the CS system.
            The pattern includes scaling, in which the length of edges in the $(n+1)$-th ring is that in the $n$-th one scaled by a factor $s$.
            Additionally, a virtual auxiliary fold parameterized by $\mu$ is inserted to introduce the constant drift term, corresponding to the orange shaded region.
            The gray shaded region represents the domain of the module.}
        \end{figure}

        In the context of two-dimensional dynamical systems, a dissipative system with a constant determinant of the Jacobian matrix is known as a CS system. Recent developments in dynamical systems theory have expanded KAM theory to include CS systems \cite{CALLEJA2013978}.
        These studies show that by adjusting a constant called the drift parameter $\mathbf{\mu}$, which offsets the expansion or contraction of phase space, quasi-periodic attractors (invariant curves) with specific frequencies can be identified. In this study, to incorporate the drift term into the module formulation of the previous work, we introduce a virtual fold that explicitly produces this drift effect (see Fig.~\ref{fig:cs_pattern}).
        Using numerical calculations, we aim to confirm the existence of these quasi-periodic attractors under the modified formulation.

        Considering the similar expansion and the drift effect, the modified map $\widetilde{\mathbf{M}}$ is defined as follows.
        
        \begin{equation}
            \widetilde{\mathbf{M}} = \boldsymbol{\Lambda} \circ \widetilde{\mathbf{M}}_{m-1} \circ \dots \circ \widetilde{\mathbf{M}}_{0}.
        \end{equation}
        Each step map $\widetilde{\mathbf{M}}_i$ is composed of the three-sphere intersection map $\widetilde{\mathbf{f}}_i$ and the connectivity transformation $\widetilde{\mathbf{g}}_i$ :
        \begin{equation}
            \widetilde{\mathbf{M}}_i = \widetilde{\mathbf{g}}_i^{k_i} \circ \widetilde{\mathbf{f}}_i,
        \end{equation}
        The component maps $\tilde{\mathbf{f}}_i$ and $\tilde{\mathbf{g}}_i$ are defined to incorporate the scaling factor $s$ only at the final step ($i=m-1$):
        \begin{eqnarray}
            \widetilde{\mathbf{f}}_i([\theta, I]^T) &:=& \begin{cases}
                \mathbf{f}_i([\theta, I]^T) \quad \mathrm{if} \ \ i < m-1, \\
                \mathbf{f}([\theta, I]^T; l_{i,L}, l_{i,R}, N, s r_{i,1}, r_{i,2}, s r_{i,3}, \sigma_i) \\ \quad \mathrm{if} \ \ i = m-1,
            \end{cases} \nonumber \\
            \widetilde{\mathbf{g}}_i([\theta, I]^T) &:=& \begin{cases}
                \mathbf{g}_i([\theta, I]^T) & \mathrm{if} \ \  i < m-1, \\
                \mathbf{g}([\theta, I]^T; s r_{i,3}, s r_{i,1}, N) & \mathrm{if} \ \  i = m-1.
            \end{cases} \nonumber
        \end{eqnarray}
        
        The adjustment map $\boldsymbol{\Lambda}$ normalizes the action variable $I$ (to match the domain of the initial state) and introduces the drift term $\mu \in \mathbb{R}$:
        \begin{equation}
            \boldsymbol{\Lambda} \left( [\theta,  I]^T \right) := \left[\theta, I/s + \mu\right]^T. \label{eq:map_lambda}
        \end{equation}
        
        Since the folding maps $\tilde{\mathbf{f}}_i$ and $\tilde{\mathbf{g}}_i$ are area-preserving (symplectic), the determinant of the Jacobian matrix of the total system depends solely on $\boldsymbol{\Lambda}$.
        Therefore, the determinant is
        \begin{equation}
            \det D\widetilde{\mathbf{M}} =  \frac{1}{s}.
        \end{equation}
        Thus, this discrete dynamical system is a CS system.

\section{\label{sec:theory}THEORETICAL ANALYSIS}
    In this section, we analyze the discrete dynamical system defined in Sec.~\ref{sec:model}.
    Because the system is constructed from symplectic (area-preserving) mappings, it admits a generating function that governs the dynamics. Analyzing this generating function is essential for clarifying the geometric and physical origins of the folding motion. However, the explicit expressions for the individual component maps are algebraically complex, making direct analytical treatment of the full discrete system infeasible.
    
    To overcome this challenge, we consider a formal asymptotic limit in which the number of modules approaches infinity
    ($N \to \infty$).
    In this continuum limit, the system becomes integrable, enabling  theoretical analysis. This approach is supported by KAM theory, which indicates that the invariant curves (quasi-periodic solutions) of the integrable system persist under sufficiently small perturbations, that is, for large but finite numbers of modules.
    
    \subsection{\label{sec:integrable}Reduction to an Integrable System}
            We begin by considering the formal large-module limit in which the number of modules
            $N$ approaches infinity ($N \to \infty$).
            From Eqs.~(\ref{eq:map_f_I}) and (\ref{eq:map_g_I}), it is immediately apparent
            that the action variable is strictly conserved:
            \begin{equation}
                I_{n+1} = I_n .
            \end{equation}
            
            Next, we examine the angular update law for the map $\mathbf{f}$
            (see Eq.~(\ref{eq:map_f_theta})).
            Let $Y_f$ and $X_f$ denote the first and second arguments of the
            $\operatorname{arctan2}(Y,X)$ function, respectively.
            In the continuum limit $N \to \infty$, explicit calculation yields
            \begin{eqnarray}
               \lim_{N \to \infty} X_f
               &=&
               \left( -e_1\beta + e_2\alpha \right)\cos\theta - e_3 \sin\theta \nonumber \\
               &=&
               \sqrt{A^2+B^2}\cos\left(\theta + \xi_f \right), \\
               \lim_{N \to \infty} Y_f
               &=&
               e_3 \cos\theta + \left( -e_1\beta + e_2\alpha \right)\sin\theta \nonumber \\
               &=&
               \sqrt{A^2+B^2}\sin\left(\theta + \xi_f \right),
            \end{eqnarray}
            where
            \begin{equation}
                A(I) = -e_1\beta + e_2\alpha, \qquad
                B(I;\sigma) = e_3 , \label{eq:AB_def}
            \end{equation}
            and the \emph{frequency} $\xi_f$ is defined by
            \begin{equation}
               \cos\xi_f = \frac{A}{\sqrt{A^2+B^2}}, \qquad
               \sin\xi_f = \frac{B}{\sqrt{A^2+B^2}} .
            \end{equation}
            Equivalently, $\xi_f$ can be expressed as
            \begin{equation}
                \xi_f = \operatorname{arctan2}(B,A). \label{eq:freq}
            \end{equation}
            Therefore, the angular update law associated with the map $\mathbf{f}$ simplifies to
            \begin{equation}
                \theta^{\prime} = \theta + \xi_f \pmod{2\pi}.
            \end{equation}
            
            We now turn to the angular update law for the connectivity map $\mathbf{g}$
            (see Eq.~(\ref{eq:map_g_theta})).
            Let $X_g$ and $Y_g$ be the quantities defined in
            Eqs.~(\ref{eq:def_Xg}) and (\ref{eq:def_Yg}), respectively.
            By multiplying these expressions by $\sin(\pi/N)$ and taking the limit
            $N \to \infty$, we obtain
            \begin{eqnarray}
                \lim_{N \to \infty} \left( \sin\left(\frac{\pi}{N}\right) X_g \right)
                &=& -4I^2 T_1 \cos\theta, \\
                \lim_{N \to \infty} \left( \sin\left(\frac{\pi}{N}\right) Y_g \right)
                &=& -4I^2 T_1 \sin\theta .
            \end{eqnarray}
            Since $\operatorname{arctan2}(Y,X)$ is invariant under the multiplication of both
            arguments by a positive scalar, we have
            \begin{equation}
                \operatorname{arctan2}(Y_g,X_g)
                =
                \operatorname{arctan2}\!\left(
                    \sin\left(\frac{\pi}{N}\right) Y_g,
                    \sin\left(\frac{\pi}{N}\right) X_g
                \right).
            \end{equation}
            Consequently, the angular update law for the map $\mathbf{g}$ is
            \begin{equation}
                \theta^{*} = \theta + \pi \pmod{2\pi}.
            \end{equation}
            
            Finally, we combine the contributions from all folding steps to determine
            the total frequency $\xi(I)$.
            Since the total map $\widetilde{\mathbf{M}}$ is defined as a composition of
            $m$ folding maps $\mathbf{f}_i$, along with multiple applications of the connectivity map $\mathbf{g}$, the angular displacements generated at each
            step accumulate additively.
            Let $\xi_{f,i}(I)$ denote the phase shift induced by the $i$-th folding map
            $\mathbf{f}_i$, and let $k_i \in \{0,1\}$ indicate whether the connectivity
            map $\mathbf{g}$ is applied at the $i$-th step.
            The total frequency $\xi(I)$ over one module is then given by
            \begin{equation}
                \xi(I)
                =
                \sum_{i=0}^{m-1} \xi_{f,i}(I)
                +
                \sum_{i=0}^{m-1} k_i \pi . \label{eq:total_freq}
            \end{equation}
            Substituting this result into the angular update law, we obtain the explicit
            form of the map in the limit $N \to \infty$:
            \begin{eqnarray}
                I_{n+1} &=& \frac{1}{s} I_n + \mu, \\
                \theta_{n+1} &=& \theta_n + \xi(I_{n+1}) \pmod{2\pi}.
            \end{eqnarray}
            If $s = 1$ and $\mu = 0$, the action variable is conserved, and the map
            reduces to an integrable action--angle system consisting of uniform
            rotations on invariant curves. The phase space is foliated by these
            invariant curves, each corresponding to a constant value of the action
            variable. When the rotation number $\xi/(2\pi)$ is irrational, the
            corresponding trajectories are quasi-periodic and dense on the associated invariant curve.
    
    \subsection{\label{sec:inv}Existence of Invariant Curves}
        Consider the symplectic map $\mathbf{M}$ in the area-preserving case (i.e., without expansion or contraction), which is governed by a generating function $S(\theta_n, I_{n+1})$.
        The discrete evolution equations are obtained from the following implicit relations:
        \begin{equation}
            I_n = \frac{\partial S}{\partial \theta_n}, \quad \theta_{n+1} = \frac{\partial S}{\partial I_{n+1}}. \label{eq:canonical_relationship}
        \end{equation}
        We express the generating function using the perturbation parameter $\varepsilon \in \mathbb{R}_{>0}$ as:
        \begin{equation}
            S(\theta_n, I_{n+1}) = \theta_n I_{n+1} + S_0(I_{n+1}) + \varepsilon S_1(\theta_n, I_{n+1}; \varepsilon),
        \end{equation}
        where $S_0(I)$ determines the integrable dynamics and $S_1$ represents the perturbation.
        In the limit $\varepsilon \to 0$, the perturbation disappears, and the system becomes an integrable map:
        \begin{equation}
            I_{n+1} = I_n, \quad \theta_{n+1} = \theta_n + \frac{d S_0}{d I_{n+1}} \label{eq:gf_corresponds_xi}.
        \end{equation}
        Here, the derivative $\xi(I) = d S_0 / d I$ corresponds to the frequency of the unperturbed system.

        We next give a formal asymptotic justification for treating the inverse module number as the perturbation parameter, i.e., $\varepsilon \sim N^{-1}$. 
        The  justification follows directly from the analytical structure of the mapping equations. As indicated in Eqs.~(\ref{eq:map_f_theta})--(\ref{eq:def_Xg}), the dependence on the module number $N$ appears exclusively through trigonometric arguments of the form $\pi/N$ or through rotation matrices $\mathbf{R}_X(2\pi/N)$ about the $X$-axis. 
        By defining $\varepsilon \equiv 1/N$, these functions become analytic with respect to $\varepsilon$ around $\varepsilon = 0$, which allows for a Taylor series expansion. 
        For example, the rotation matrix can be expanded as
        \begin{equation}
            \mathbf{R}_X(2\pi \varepsilon) = \mathbf{E} + 2\pi \varepsilon \mathbf{J}_X + O(\varepsilon^2),
        \end{equation}
        where $\mathbf{E}$ denotes the $3\times 3$ identity matrix, and $\mathbf{J}_X$ is the generator of infinitesimal rotations about the $X$-axis, given explicitly by
        \begin{equation}
            \mathbf{J}_X = 
            \begin{pmatrix}
                0 & 0 & 0 \\
                0 & 0 & -1 \\
                0 & 1 & 0
            \end{pmatrix}.
        \end{equation}
        Here, the first-order term naturally scales with $N^{-1}$.
        Furthermore, while the map contains terms that appear singular, such as $\cot(\pi/N)$, these are regularized by the normalization involving $\sin(\pi/N)$. 
        For instance, the identity
        \begin{equation}
            \sin(\pi \varepsilon) \cot(\pi \varepsilon) = \cos(\pi \varepsilon) = 1 - \frac{\pi^2}{2}\varepsilon^2 + O(\varepsilon^4)
        \end{equation}
        demonstrates that such singularities are eliminated, resulting in expressions that are analytic in $\varepsilon$.
        Consequently, the total map $\mathbf{M}$ naturally allows a perturbation expansion of the form
        \begin{equation}
            \mathbf{M} = \mathbf{M}_0 + \varepsilon \mathbf{M}_1 + O(\varepsilon^2),
        \end{equation}
        where the zeroth-order term $\mathbf{M}_0$ exactly recovers the integrable limit derived in Sec.~\ref{sec:integrable}. 
        This analysis confirms that the inverse module number $\varepsilon = 1/N$ functions as a valid perturbation parameter for the subsequent KAM analysis.
        
        To apply the KAM theorem, we additionally assume that the \emph{nondegeneracy condition} \begin{equation} 
            \frac{d\xi}{dI} \neq 0
            \end{equation}
        is satisfied. It is important to note that the twist condition is a local requirement with respect to the action variable $I$ and does not assume that $\xi(I)$ is monotonic over the entire phase space. In particular, even if $\xi(I)$ possesses extrema at certain values of $I$, KAM theory is applicable at other values of $I$ where the nondegeneracy condition $d\xi/dI \neq 0$ is satisfied. We also assume that $\xi(I_0)$ satisfies a \emph{Diophantine condition}.
        Since the angular variable $\theta$ is defined on the circle $\mathbb{T}^1 = \mathbb{R}/2\pi\mathbb{Z}$, the non-resonance condition must be evaluated modulo $2\pi$.
        Specifically, we require that there exist constants $C > 0$ and $\tau>1$ such that
        \begin{equation}
            |k\,\xi(I_0) - 2\pi l| \ge \frac{C}{|k|^\tau}, \quad \forall k\in\mathbb{Z}\setminus\{0\}, \ \forall l \in \mathbb{Z}.
        \end{equation}
        Under these assumptions, the KAM theorem guarantees that, for sufficiently small perturbations, the invariant curve corresponding to $I=I_0$ persists, carrying quasi-periodic dynamics with frequency $\xi(I_0)$.
        
        Conversely, systems in which the frequency profile $\xi(I)$ exhibits extrema are known as \emph{nontwist systems}. 
        In these systems, the monotonicity of the frequency is disrupted, resulting in different action values on either side of the extremum corresponding to the same frequency.  
        This degeneracy leads to dynamical phenomena distinct from those observed in standard twist maps, such as \emph{reconnection} and bifurcation of invariant curves near the extremum, often referred to as the \emph{shearless curve}\cite{DELCASTILLONEGRETE19961}. 
        In the map investigated in this study, the violation of the twist condition can be readily controlled.
        
        As indicated in Eq.~(\ref{eq:freq}), the frequency profile $\xi(I)$ depends on the parameter $\sigma$ within the term $B(I;\sigma)$. 
        By adjusting $\sigma$, one can superpose terms with competing signs, thereby inducing local maxima or minima in $\xi(I)$. 
        As a result, this process transforms the system into a nontwist map.
        
        In these discrete dynamical systems, stable elliptical islands generally form at resonances, particularly where the frequency approaches zero ($\xi(I) \approx 0$). Previous studies examined configurations with at most one stability island, constrained by a frequency profile intersecting the resonance condition only once. We show that modifying the origami parameters enables a more intricate frequency topology, with zeros at multiple action values.
    
       A critical challenge in these systems is that, even in nontwist scenarios, the resonance condition $\xi(I) = 0$ may be disrupted by perturbations; the frequency curve can shift away from the zero axis or enter regions where no physically feasible solutions exist. In contrast, by leveraging the nontwist characteristics of the current system and precisely adjusting the positions of frequency extrema, we can ensure multiple stable intersections with the zero-frequency axis. This multiplicity of resonances results in the coexistence of several distinct elliptical islands within the phase space, offering a level of design flexibility not previously observed in tubular origami tessellations.

        The classical KAM theory guarantees the persistence of invariant curves only in conservative (Hamiltonian) systems. 
        The generalized origami map $\widetilde{\mathbf{M}}$ with similarity expansion (Sec.~\ref{sec:cs}) is dissipative, with constant Jacobian $1/s \neq 1$, so that invariant tori generally break down, leaving isolated periodic or chaotic attractors.
        
        Recent extensions of KAM theory to CS systems show that quasi-periodic invariant tori can persist if a drift parameter is adjusted to compensate for contraction or expansion. 
        In our tubular origami model, the auxiliary fold parameter $\mu$ in Eq.~(\ref{eq:map_lambda}) serves as this drift. 
        Thus, for sufficiently small perturbations (large $N$), stable quasi-periodic attractors are expected for appropriate values of $\mu$, providing a theoretical basis for the numerical results presented in the next section.

    \subsection{\label{sec:genfun}The Generating Function and Variational Structure}
        The dynamical system governed by the map $\mathbf{M}$ (without expansion) is conservative. For symplectic mappings, a generating function exists locally. Although obtaining an explicit expression in the general case is analytically infeasible due to the map's algebraic complexity, we have identified the generating function in the integrable limit. Importantly, this function can be interpreted as the total discrete mean curvature. In this section, we present the derivation and analysis of this variational structure.
    
        Since the total frequency $\xi(I)$ is constructed as a linear sum of the phase shifts contributed by individual zigzags (as shown in Eq.~(\ref{eq:total_freq})), the total generating function $S_0(I)$ naturally decomposes into a sum of local generating functions associated with each folding step.
        Therefore, it is sufficient to focus our analysis on the generating function $S_f(I)$ corresponding to a single folding map $\mathbf{f}$.
        Using the definition $\xi_f = \operatorname{arctan2}(B, A)$, the function $S_f(I)$ is derived via integration by parts:
        \begin{eqnarray}
            S_f(I) &=& \int \xi_f(I) \, dI = I \xi_f - \int I \frac{d\xi_f}{dI} \, dI. \label{eq:S_f_first_def}
        \end{eqnarray}
        We define the integrand of the second term as $\psi(I) := I (d\xi_f/dI)$.
        Using the explicit form of the derivative, $\psi(I)$ can be written in terms of a rational function $R(I)$ as:
        \begin{equation} 
            \psi(I) = \frac{I(B^{\prime}A-A^{\prime}B)}{A^2+B^2} = \frac{\sigma}{\beta\sqrt{\gamma}} R(I). \label{eq:def_psi}
        \end{equation}
        (See Appendix~\ref{app:derivation} for the detailed algebraic derivation).
        
        The integrability of the term $\psi(I)$ is strictly governed by the algebraic structure of the denominator, specifically the irrational factor $\beta\sqrt{\gamma}$.
        Let $P(I)$ denote the polynomial inside the square root appearing in this factor.
        According to the classical theory of integration\cite{silverman2009arithmetic}, the nature of the primitive function is determined by the degree of $P(I)$, denoted by $d = \deg(P)$.
        If $d \le 2$, the integral can be evaluated explicitly using elementary functions, typically involving inverse trigonometric or logarithmic functions. However, if $d = 3$ or $d = 4$, the integral generally cannot be expressed in terms of elementary functions and are classified as \emph{elliptic integrals}.
        For $d \ge 5$, the problem involves \emph{hyperelliptic integrals}, which represent a higher-complexity class of transcendental functions.
    
        In the present origami system, the degree $d$ is directly dictated by the geometric symmetry of the module.
        As detailed in Appendix~\ref{app:polynomial_degree}, an analysis of the explicit form of $P(I)$ reveals that the integrability class bifurcates depending on the relationship between the left and right crease lengths.
        In the general asymmetric case ($l_L \neq l_R$), the polynomial is found to be of degree $d=4$, implying that the generating function is intrinsically non-elementary and is described by elliptic integrals.
        In contrast, when the module possesses symmetry ($l_L = l_R$), the algebraic structure simplifies significantly; the coefficients of the higher-order terms vanish, reducing the polynomial to degree $d=2$.
        Consequently, in this symmetric case, the variational structure can be expressed entirely in terms of elementary functions.
        In this study, we specifically focus on this symmetric configuration to derive an explicit analytical representation of the generating function.

        Relegating the detailed integration steps to Appendix~\ref{app:integration_result}, the integral of $\psi(I)$ for the symmetric case (where $l_L = l_R = l$ and $r_1 = r_3$) yields the following explicit closed-form expression:
        \begin{equation}
            -\int\psi(I)dI = l\Phi_l + r_1\Phi_{r_1} + r_2\frac{\Phi_{r_2}}{2} + \mathrm{Const.}, \label{eq:psi_result}
        \end{equation}
        where $\Phi_l, \Phi_{r_1},$ and $\Phi_{r_2}$ represent the exterior dihedral angles associated with the crease lengths $l, r_1,$ and $r_2$, respectively.
        
        Finally, substituting Eq.~(\ref{eq:psi_result}) back into the expression derived via integration by parts, we obtain the total generating function $S_f(I)$:
        \begin{equation}
            S_f(I) = \frac{1}{2} \left(2I\xi_f + 2l\Phi_l + 2r_1\Phi_{r_1} + r_2\Phi_{r_2}\right) + \mathrm{Const.}
        \end{equation}
        In this expression, the first term $I\xi_f$ represents half the product of the crease length $2I$ and its associated exterior dihedral angle $\xi_f$ (see Appendix~\ref{app:geometric_interpretation}).
        Consequently, the generating function $S_f(I)$ corresponds to half the sum of the products of the lengths and the exterior dihedral angles of all edges forming the tetrahedron $U_1U_2U_3-V_\sigma$ (as shown in Fig.~\ref{fig:parametrization_zigzag}).
        From the perspective of discrete differential geometry, this quantity is identified as the total discrete mean curvature $H$ of the discrete surface, defined as follows\cite{Bobenko2008Book}:
        \begin{equation}
            H = \frac{1}{2} \sum_{j \in \mathcal{E}} l_j \Phi_j, \label{eq:mean_curvature}
        \end{equation}
        where the summation extends over all edges $j \in \mathcal{E}$ of the polyhedron, with $l_j$ denoting the edge length and $\Phi_j$ representing the corresponding exterior dihedral angle. This formula is the standard definition of total mean curvature for polyhedral surfaces in discrete differential geometry. Notably, this discrete quantity corresponds to the integral of mean curvature over a smooth surface, denoted as $\int H \, dA$.
        This relationship establishes a geometric law: the conserved quantity governing the integrable dynamics of tubular origami is the total discrete mean curvature.

        Synthesizing the results derived above, we extend this local formulation to the entire module structure.
        The total generating function $S_0(I)$ is obtained by summing the contributions from both the folding maps $\mathbf{f}_i$ and the re-parametrization maps $\mathbf{g}_i$.
        Recalling that the map $\mathbf{g}$ acts as a pure phase shift of $\pi$ in the continuum limit (i.e., $\xi_g = \pi$), its generating function is given by the linear term $S_{g,i}(I) = k_i \pi I$, where $k_i \in \{0, 1\}$ is the connectivity index.
        Consequently, the total generating function takes the form:
        \begin{equation}
            S_0(I) = \sum_{i=0}^{m-1} S_{f,i}(I) + \sum_{i=0}^{m-1} k_i \pi I = H_{\mathrm{module}}(I)
        \end{equation}
        where $H_{\mathrm{module}}(I)$ is the total discrete mean curvature in one module.
               To formulate the variational principle for the integrable dynamics, it is more natural to use a generating function depending on two consecutive angle variables rather than on two action variables.
        Indeed, in the integrable limit, different values of the action variable correspond to different invariant curves, and therefore a sequence described by transitions between different action values does not represent a single orbit of the integrable system.
        We therefore introduce the type-I generating function, or discrete Lagrangian, by the Legendre transform of $S(\theta_n,I_{n+1})$:
        \begin{equation}
            L(\theta_n,\theta_{n+1})
            =
            S(\theta_n,I_{n+1}) - I_{n+1}\theta_{n+1},
            \label{eq:type1_generating_function}
        \end{equation}
         where $I_{n+1}=I_{n+1}(\theta_n,\theta_{n+1})$ is understood as the action variable determined through the Legendre transform.
        This transformation is taken locally on intervals where the twist condition $d\xi/dI \neq 0$ holds, so that the relation between the angular increment and the action variable can be locally inverted.
        With this convention, the type-I generating function becomes
        \begin{equation}
            L(\theta_n,\theta_{n+1})
            =
            S_0(I_{n+1}) - I_{n+1}\xi(I_{n+1}).
            \label{eq:lagrangian_legendre}
        \end{equation}
        Since $S_0(I)$ is identified with the total discrete mean curvature $H_{\mathrm{module}}(I)$, the Legendre transform subtracts the contribution associated with the edge corresponding to the angular increment.
        As a result, the discrete Lagrangian is expressed only in terms of the remaining discrete mean-curvature contributions of the module.

        Using this type-I generating function, we define the discrete action functional for a sequence of angle variables $\{\theta_n\}_{n=0}^{M}$ as
        \begin{equation}
            \mathcal{A}
            =
            \sum_{n=0}^{M-1}
            L(\theta_n,\theta_{n+1}).
            \label{eq:discrete_action_theta}
        \end{equation}
        
        The trajectory of the folding dynamics is determined by the principle of stationary action,
        $\delta \mathcal{A}=0$, subject to fixed endpoint conditions for the angle variables.
        The corresponding discrete Euler--Lagrange equation is
        \begin{equation}
            \frac{\partial L(\theta_{n-1},\theta_n)}{\partial \theta_n}
            +
            \frac{\partial L(\theta_n,\theta_{n+1})}{\partial \theta_n}
            =
            0.
            \label{eq:discrete_EL_theta}
        \end{equation}
        In area-preserving variational systems, local matching conditions often
        appear as geometric rules; a familiar example is the reflection law in
        billiard systems\cite{Tabachnikov2005Billiards}.
        Similarly, the variational equation obtained here can be viewed as a
        local compatibility condition between adjacent origami modules, rather
        than as a condition imposing constant mean curvature over the entire
        folded surface.
        \begin{figure}
            \includegraphics[width=\columnwidth]{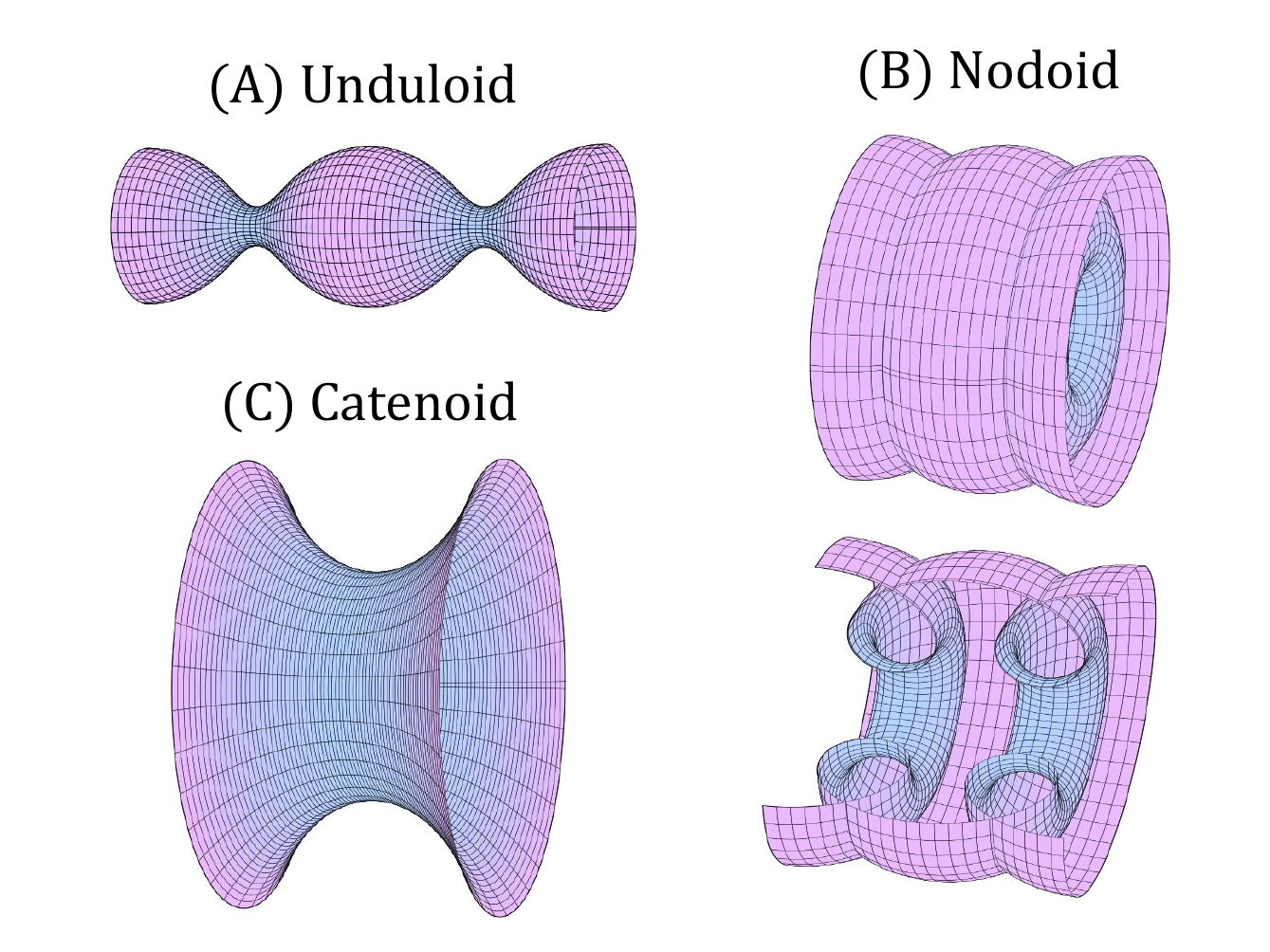}
            \caption{\label{fig:delaunay_surfaces}
                Representative classifications of constant mean curvature (CMC) surfaces. These surfaces are obtained as solutions to the variational problem of maximizing volume subject to a constant surface area. When rotational symmetry is imposed, these CMC surfaces are known as Delaunay surfaces. Examples include:
                (A) Unduloid, a smooth continuous surface;
                (B) Nodoid, which exhibits self-intersections, along with its cross-section;
                (C) Catenoid, the limiting shape where the mean curvature becomes zero.
            }
        \end{figure}
        In the theory of continuous surfaces, constant mean curvature (CMC)
        surfaces\cite{Kenmotsu2003SurfacesCMC}, as illustrated in
        Fig.~\ref{fig:delaunay_surfaces}, form a well-known class of surfaces arising
        from variational principles involving surface area and enclosed volume.
        The folded origami configurations considered in this study show geometric
        similarities to such CMC surfaces.
        This observation suggests that the generating-function formulation developed
        here may be understood as a discrete, curvature-based variational analogue of
        CMC-type geometry in tubular origami dynamics.
        
\section{\label{sec:numerical}NUMERICAL EXPERIMENTS}
    Based on the theoretical analysis in Sec.~\ref{sec:theory}, we conduct numerical experiments to verify the dynamical behavior of the system. We aim to confirm the existence of invariant curves (KAM tori) for sufficiently large $N$ in both symplectic and CS regimes. To enhance understanding of the dynamics, we visualize the corresponding 3D origami structures for representative trajectories in phase space. During numerical iterations, there are instances where no solution satisfies the geometric constraint; specifically, when the intersection of three spheres at a crease vertex yields no solution. We refer to these cases as \emph{finite solutions}. In phase diagrams, any trajectory that reaches a finite solution is terminated and therefore omitted from the plot. In addition, we treat as finite solutions those trajectories that pass through the singular point $I = 0$, namely, cases in which the vertices $U_1$ and $U_3$ interchange.
    
    \subsection{\label{sec:kam}KAM Tori in Symplectic Dynamical Systems}
    
    We performed numerical experiments on symplectic dynamical systems governed by the mappings
    $\mathbf{M} = \mathbf{f}_3 \circ \mathbf{f}_2 \circ \mathbf{f}_1$ and
    $\mathbf{M} = \mathbf{f}_3 \circ (\mathbf{g}_2 \circ \mathbf{f}_2) \circ \mathbf{f}_1 \circ (\mathbf{g}_0 \circ \mathbf{f}_0)$,
    which correspond to the origami modules shown in Fig.~\ref{fig:examples_of_modules}.
    As discussed in Sec.~\ref{sec:inv}, by appropriately tuning the geometric parameters, these systems can be designed to realize both
    twist systems, in which the frequency function $\xi(I)$ is monotonic, and
    nontwist systems, in which $\xi(I)$ exhibits local extrema.
    Although candidate parameter sets can be identified numerically from the profile of $\xi(I)$,
    inappropriate choices of edge lengths may result in self-intersections or geometric collapse. Therefore, throughout this study, we selected parameters that ensure physically admissible, non-self-intersecting origami configurations near elliptic fixed points.
    
    \begin{figure}
        \includegraphics[width=\columnwidth]{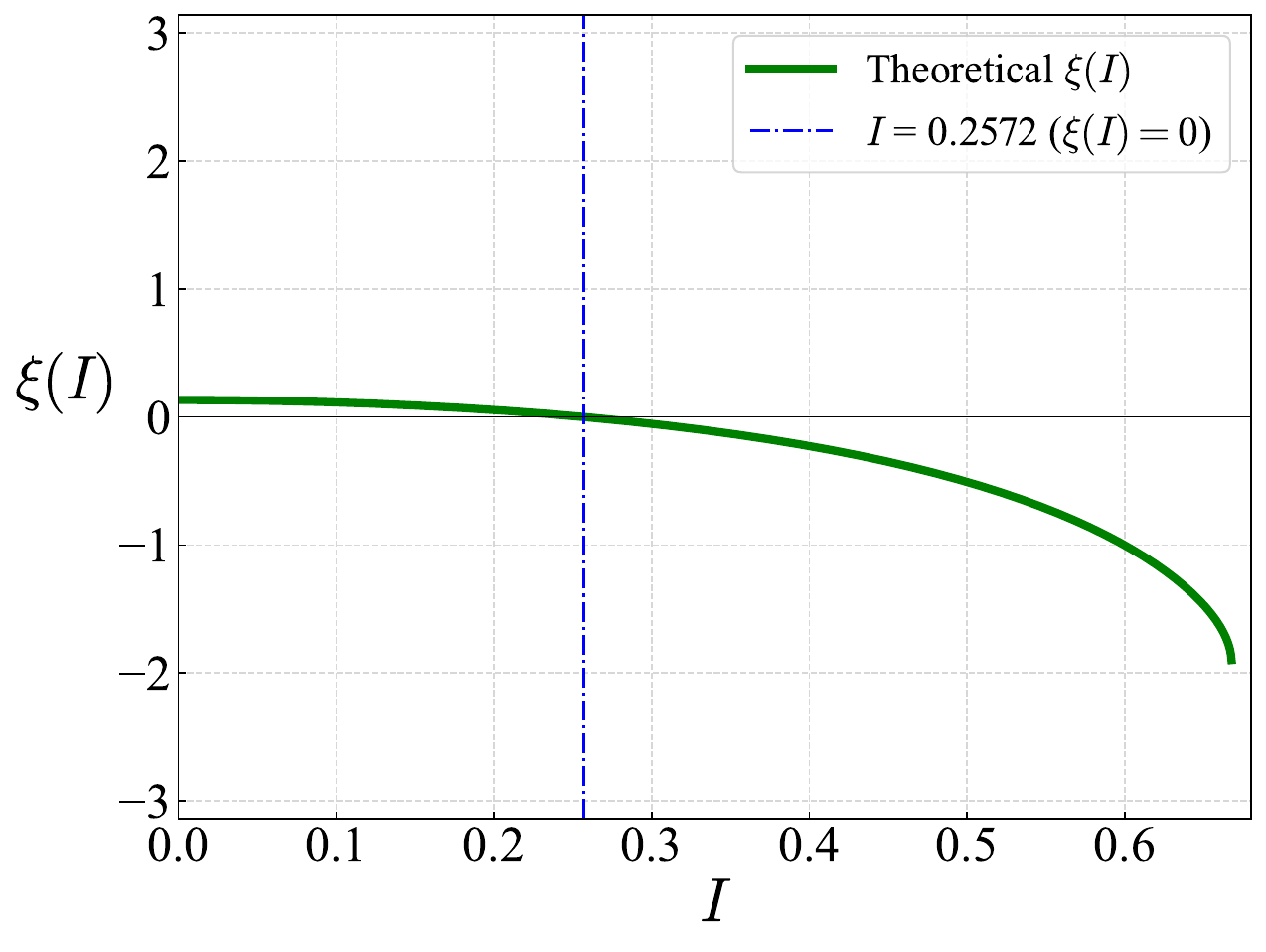}
        \caption{\label{fig:origami_xi_profile_1}
            Numerical profile of the frequency function $\xi(I)$ for the basic module.
            The solid curve represents the theoretical value of $\xi(I)$.
            The function is strictly monotonic ($d\xi/dI < 0$) over the entire valid domain, indicating that the twist condition is satisfied.
            The blue dash-dotted line marks the location where $\xi(I)=0$, at $I \approx 0.2572$.
        }
    \end{figure}
    
\begin{figure*}
        \includegraphics[width=16.5cm]{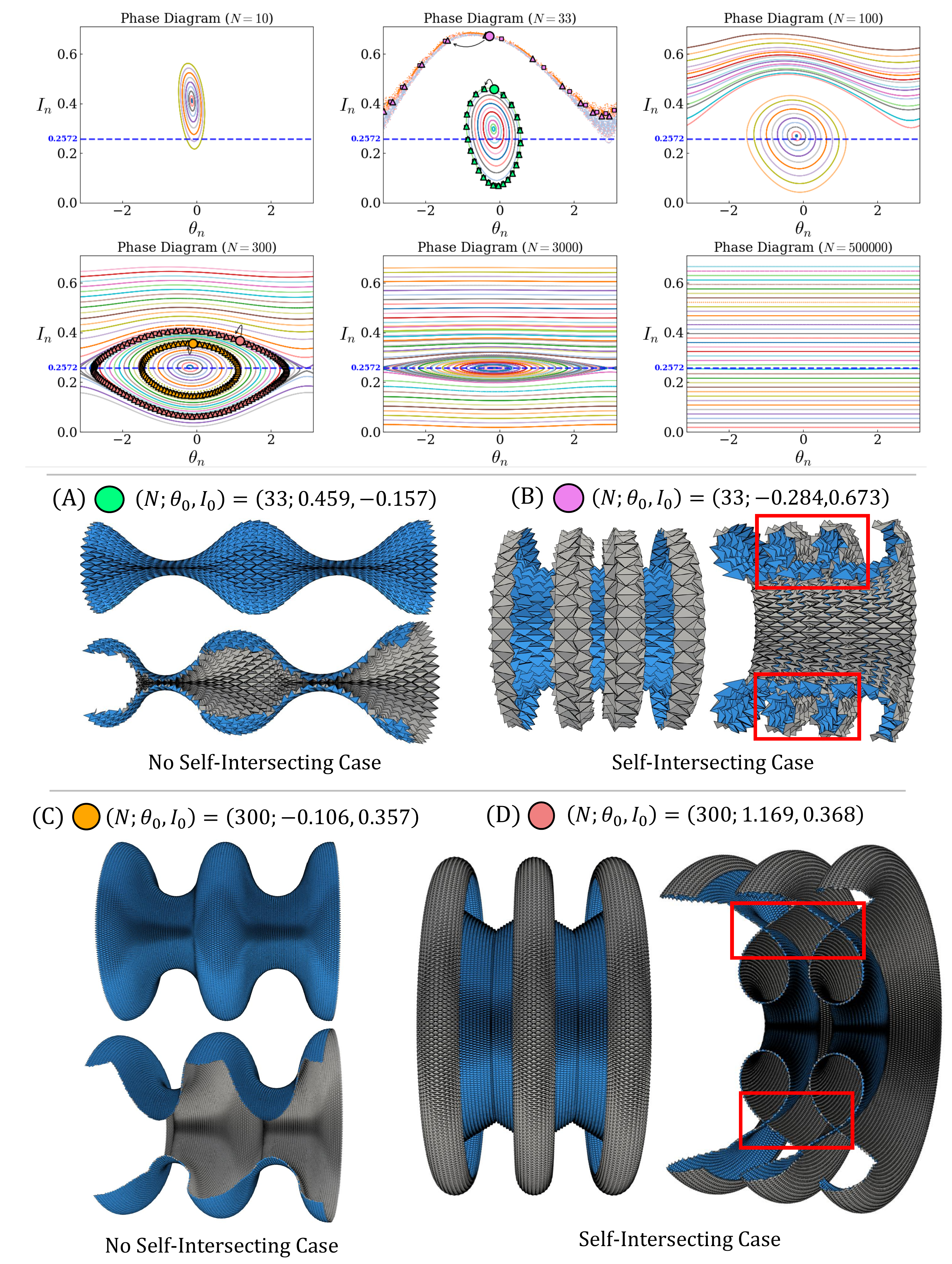}
        \caption{\label{fig:numerical_result_1}
            Phase diagrams of the mapping $\mathbf{M}=\mathbf{f}_2 \circ \mathbf{f}_1 \circ \mathbf{f}_0$, computed for module numbers
            $N = 10, 33, 100, 300, 3000$, and $500000$.
            Orbits corresponding to different initial conditions are plotted in different colors.
            For selected trajectories, the initial point $(\theta_0, I_0)$ is indicated by a circle marker ($\circ$),
            while points belonging to the first and second periods are marked by triangles ($\triangle$) and squares ($\square$), respectively.
            The lower panels show the corresponding three-dimensional reconstructed origami configurations.
            The geometric parameters are set to $(\sigma_0, \sigma_1, \sigma_2) = (M, M, M)$,
            and the edge lengths $(l_{i,L}, l_{i,M}, l_{i,R})$ for steps $i=0,1,2$ are fixed at
            $(0.75, 0.69, 0.75)$, $(1.11, 0.58, 1.11)$, and $(0.90, 1.49, 0.90)$, respectively.
        }
    \end{figure*}

    We first present numerical results for the twist system characterized by a monotonic frequency, as shown in Fig.~\ref{fig:origami_xi_profile_1}.
    The corresponding phase diagrams are displayed in Fig.~\ref{fig:numerical_result_1}.
    The numerical results suggest that, as the number of modules $N$ increases, large regions of the plotted phase space are occupied by apparently invariant curves. For sufficiently large $N$, these invariant curves remain, although they are slightly deformed from the integrable limit. This behavior is consistent with the KAM-type persistence of nonresonant
    invariant curves under small perturbations.
    
    A closer inspection of Fig.~\ref{fig:numerical_result_1} reveals a strong correspondence between the phase-space structure
    and the theoretical frequency profile.
    For large module numbers, such as $N=300$ and $3000$, the center of the primary elliptic island is located at
    $I \approx 0.2572$, in excellent agreement with the zero-frequency condition $\xi(I)=0$
    predicted by the theoretical curve in Fig.~\ref{fig:origami_xi_profile_1}.
    For smaller values of $N$, where the discrete system experiences stronger perturbative effects, the fixed-point shifts slightly away from this theoretical value; nevertheless, it remains near the predicted location. More fundamentally, in this origami-based system, physically realizable foldable regions correspond predominantly to neighborhoods of vanishing frequency. This observation provides a clear design principle: by tuning the geometric parameters to control the profile of $\xi(I)$,
    one can effectively engineer the location of stable, foldable domains in phase space.
    
    We further examine origami configurations associated with several representative trajectories.
    Trajectories lying on elliptic quasi-periodic orbits, shown in Fig.~\ref{fig:numerical_result_1}(A) and (C),
    do not exhibit geometric failures such as self-intersections and are therefore physically foldable in real space.
    However, we find that even elliptic quasi-periodic orbits that are commonly regarded as belonging to foldable regions
    may give rise to self-intersections, as illustrated in Fig.~\ref{fig:numerical_result_1}(D).
    In contrast, the origami configuration corresponding to an invariant curve
    (Fig.~\ref{fig:numerical_result_1}(B)) clearly exhibits self-intersection.
    Although this trajectory appears chaotic due to perturbative effects,
    it in fact retains quasi-periodic motion.
    For stronger perturbations, such as in the case $N=10$,
    the invariant curve disappears altogether,
    indicating that the perturbed trajectory is driven into a region where
    physically admissible geometric solutions no longer exist.
    
    \begin{figure}
        \includegraphics[width=\columnwidth]{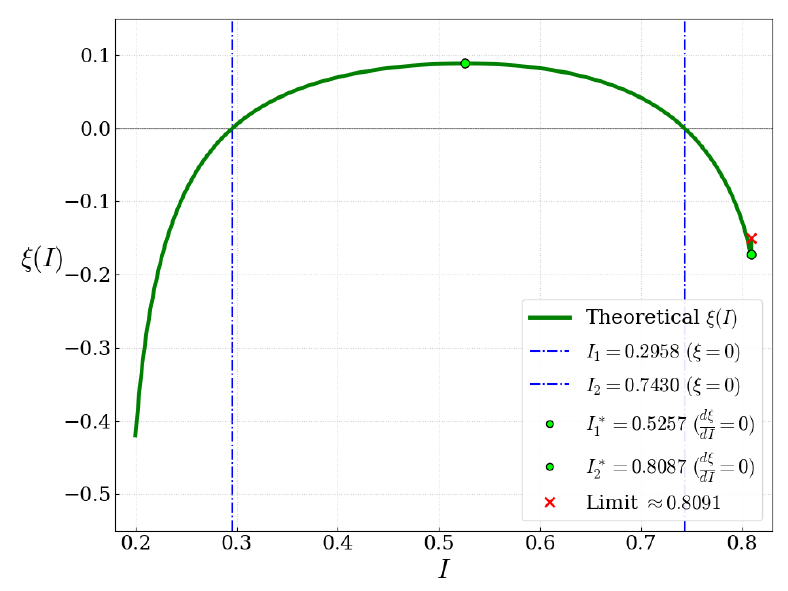}
        \caption{\label{fig:origami_xi_profile_2}
				Theoretical profile of the frequency function $\xi(I)$ for the nontwist case. In contrast to the twist system, $\xi(I)$ is non-monotonic and exhibits local extrema (shearless points) at $I_1^* \approx 0.5257$ and $I_2^* \approx 0.8087$, indicating a violation of the global twist condition. The profile clearly intersects the zero-frequency axis at two actions, $I_1 \approx 0.2958$ and $I_2 \approx 0.7430$. The curve terminates at the geometric limit $I_{\text{limit}} \approx 0.8091$ just after the second shearless point. Although a third intersection is absent in this theoretical limit, the proximity of the curve to $\xi=0$ near the boundary allows for the emergence of a third island structure under finite-$N$ perturbations.
        }
    \end{figure}
    
\begin{figure*}
        \includegraphics[width=16.5cm]{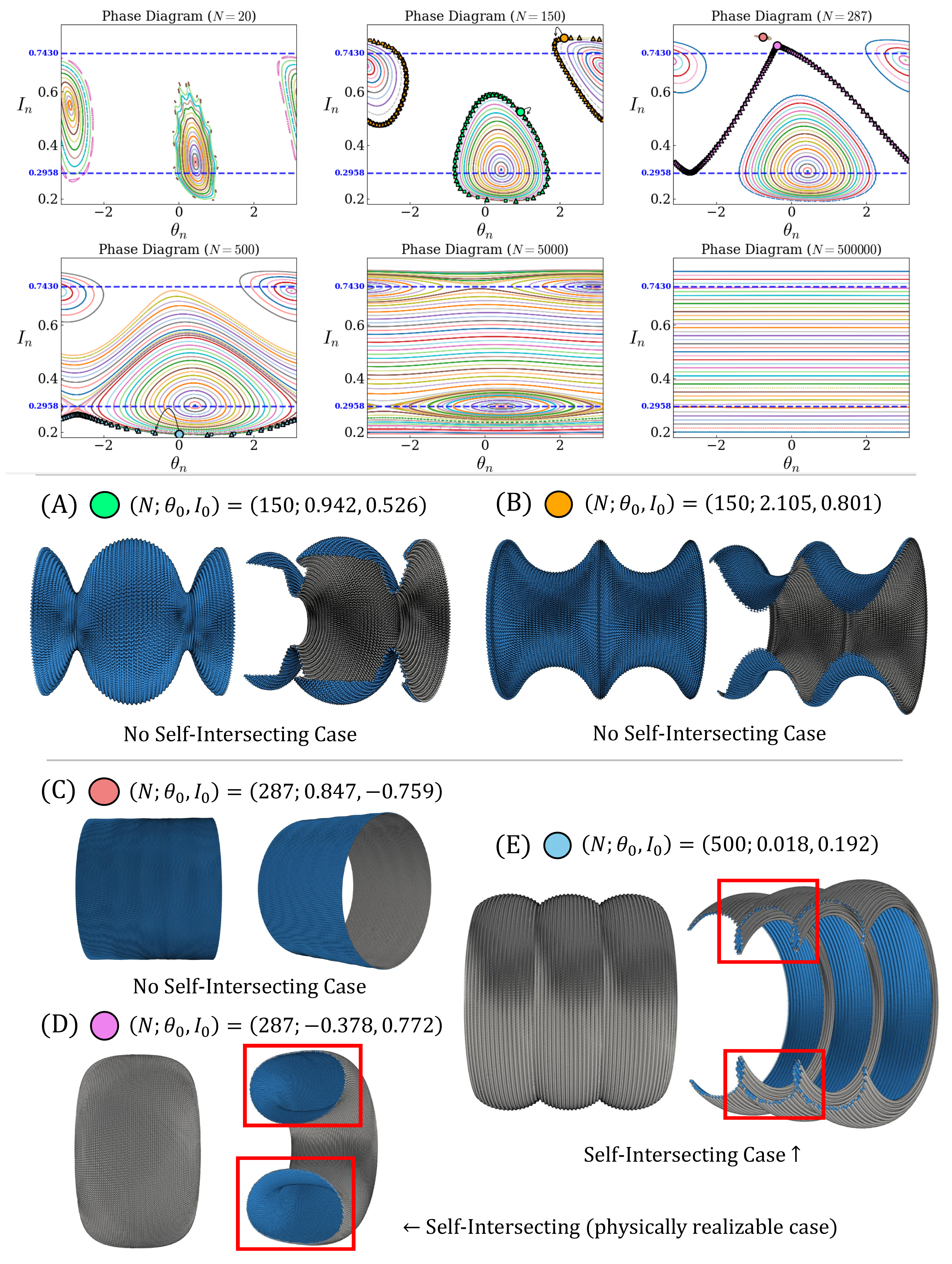}
        \caption{\label{fig:numerical_result_2}
            Phase diagrams of the mapping
            $\mathbf{M} = \mathbf{f}_3 \circ (\mathbf{g}_2 \circ \mathbf{f}_2) \circ \mathbf{f}_1 \circ (\mathbf{g}_0 \circ \mathbf{f}_0)$,
            computed for module numbers $N = 20, 150, 287, 500, 5000$, and $500000$.
            Orbits corresponding to different initial conditions are plotted in different colors.
            The initial point $(\theta_0, I_0)$ is indicated by a circle marker ($\circ$),
            while points belonging to the first and second periods are denoted by triangles ($\triangle$)
            and squares ($\square$), respectively.
            The lower panels show the corresponding three-dimensional reconstructed origami configurations.
            The geometric parameters are set to $(\sigma_0, \sigma_1, \sigma_2, \sigma_3) = (V, V, M, M)$,
            and the edge lengths $(l_{i,L}, l_{i,M}, l_{i,R})$ for steps $i=0,1,2,3$ are fixed at
            $(0.8552, 0.9913, 0.9748)$, $(1.1342, 1.1798, 1.0774)$,
            $(1.1932, 0.9922, 0.8522)$, and $(0.9105, 0.9236, 1.0664)$, respectively.
        }
    \end{figure*}

    Next, we consider the nontwist case in which $\xi(I)$ exhibits local extrema, as shown in Fig.~\ref{fig:origami_xi_profile_2}. The corresponding numerical results are presented in Fig.~\ref{fig:numerical_result_2}. As in the twist case, the phase space becomes increasingly populated with invariant curves as $N$ increases. In this specific configuration, the theoretical frequency profile $\xi(I)$ admits two roots, $I_1 \approx 0.2958$ and $I_2 \approx 0.7430$, where $\xi(I)=0$, corresponding to the two primary elliptic islands observed in the phase space. Notably, the profile terminates at the geometric limit $I_{\text{limit}} \approx 0.8091$ immediately after the second shearless point. While the theoretical curve does not strictly intersect the zero-axis a third time, it remains in extreme proximity to $\xi=0$ at this boundary. As a result, the topological instability associated with the reconnection phenomenon near the shearless region, coupled with finite-$N$ perturbations, allows for the emergence of a third, albeit marginal, stable region at the edge of the phase space. 
    
    Furthermore, in Numerical Experiment 2, where the map is constructed from a longer sequence of compositions, the impact of non-linearity is significantly amplified.
    With the increased complexity of the mapping, the cumulative perturbative effects become much stronger than in the simpler cases.
    As a result, at specific $I$ values of the zero profile the center of the stable region deviates significantly from the theoretical prediction.
    Indeed, numerical results confirm that the minute elliptic island observed in Fig.~\ref{fig:numerical_result_2}(C) corresponds to this value $I_3$; it is barely discernible in the phase diagram as a result of these combined effects.


     \begin{figure*}
                \includegraphics[width=16.5cm]{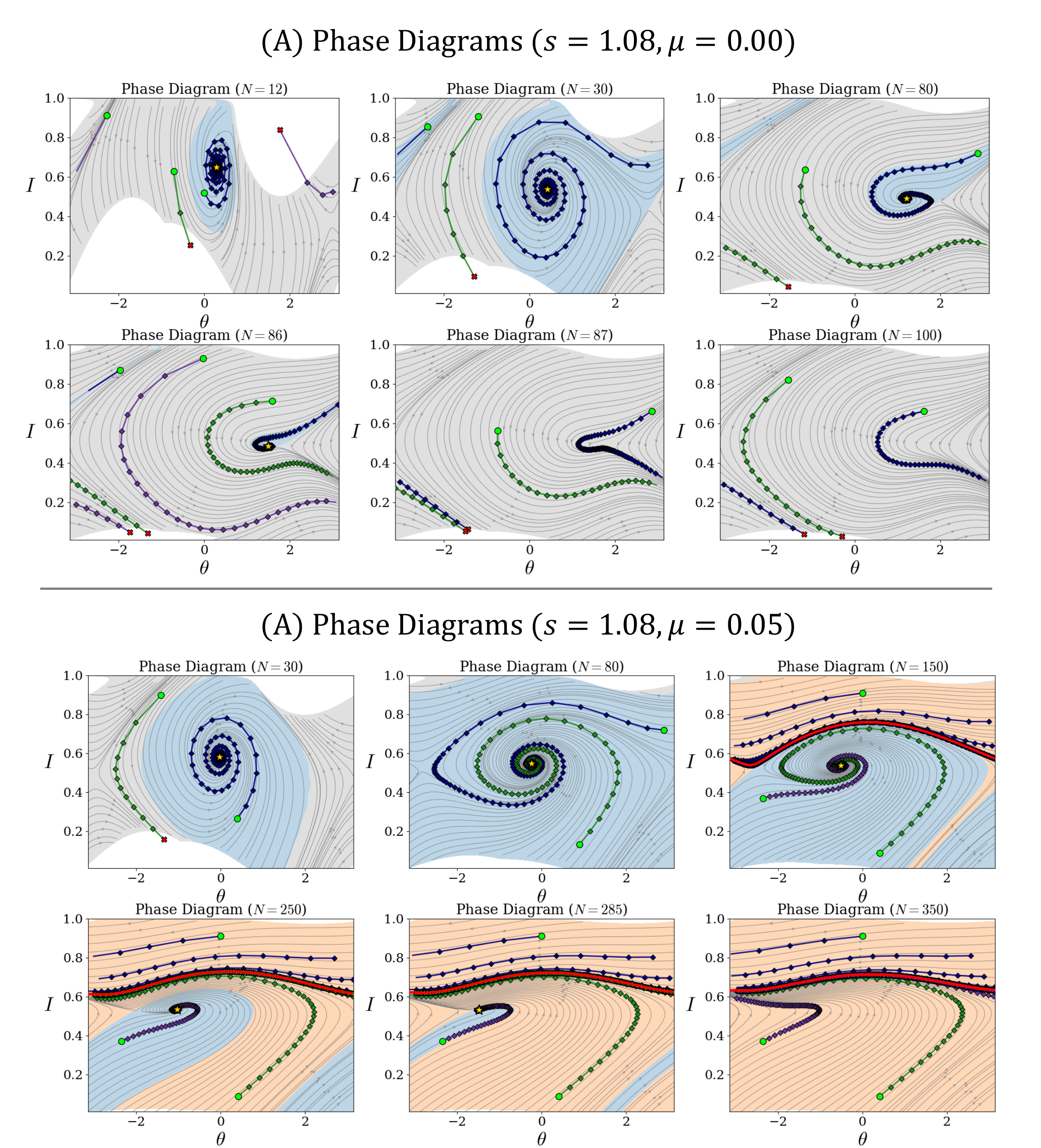}
                \caption{\label{fig:numerical_result_3}
                    Phase diagrams of the conformally symplectic mapping $\widetilde{\mathbf{M}} =\boldsymbol{\Lambda} \circ \mathbf{f}_2 \circ \mathbf{f}_1 \circ \mathbf{f}_0$, computed for various module numbers $N$.
                    Different orbits are plotted in distinct colors.
                    For each trajectory, the initial point $(\theta_0, I_0)$ is marked by a green circle ($\circ$), and subsequent iterations are denoted by diamonds ($\diamond$).
                    The system exhibits two distinct types of attractors: a point attractor, indicated by a star marker ($\star$) with its basin of attraction shaded in light blue, and a quasi-periodic attractor, represented by a solid red line with its basin shaded in orange.
                    Background streamlines visualize the flow direction of the dynamics.
                    The gray regions correspond to domains where trajectories eventually terminate (finite solutions) due to geometric constraints, while the white regions represent areas where the map is immediately undefined.
                    The geometric parameters are set to $(\sigma_0, \sigma_1, \sigma_2) = (M, M, M)$, and the edge lengths $(l_{i,L}, l_{i,M}, l_{i,R})$ for steps $i=0, 1, 2$ are fixed at $(1.0, 0.91, 1.0)$, $(\sqrt{2}, 0.91, \sqrt{2})$, and $(\sqrt{2}, 1.91, \sqrt{2})$, respectively.
                    }
        \end{figure*}

    We also examine the origami configurations associated with the invariant curves shown in
    Figs.~\ref{fig:numerical_result_2}(D) and (E).
    As in the twist regime, such trajectories generally correspond to configurations exhibiting self-intersections.
    However, the case shown in Fig.~\ref{fig:numerical_result_2}(D) constitutes a notable exception.
    Although the global trajectory involves self-intersection,
    we find that over a single period this solution yields a physically constructible toroidal shape
    without self-intersection.
    From a dynamic perspective, this orbit can be identified as a
    \emph{meandering invariant curve}, a structure unique to nontwist systems\cite{DELCASTILLONEGRETE19961}.
    Its undulated geometry indicates that the trajectory is strongly influenced by
    the topological constraints imposed by the pair of elliptic islands located on either side.
    
    Furthermore, Figs.~\ref{fig:numerical_result_2}(A) and (B) show origami configurations
    corresponding to orbits within this pair of elliptic islands.
    These two islands create contrasting geometric configurations
    that exhibit a striking geometric duality.
    The structure in Fig.~\ref{fig:numerical_result_2}(A) features a deeply indented, concave neck,
    resulting in a smooth, bulbous form.
    In contrast, the structure in Fig.~\ref{fig:numerical_result_2}(B) is characterized by sharp,
    outward-pointing cusps.
    The smooth indentation of the former is thus geometrically inverted
    into the singular protrusion of the latter.
    This correspondence suggests that the nontwist bifurcation effectively separates the dynamics
    into two distinct folding modes, which may be described as ``concave'' and ``convex.''

    \subsection{\label{sec:attractor}A KAM Attractor in the CS Dynamical Systems}
        Finally, we analyze the existence of quasi-periodic attractors predicted by KAM theory for CS systems.
        Numerical simulations confirm the presence of such attractors, and the results for the mapping $\widetilde{\mathbf{M}} = \boldsymbol{\Lambda} \circ \mathbf{f}_2 \circ \mathbf{f}_1 \circ \mathbf{f}_0$ are summarized in Fig.~\ref{fig:numerical_result_3}.

        Figure~\ref{fig:numerical_result_3}(A) presents the results for the case with zero drift ($\mu = 0$).
        Under this condition, we observe that the physically valid phase space is progressively annihilated as the module number $N$ increases.
        Eventually, for sufficiently large $N$, the system becomes dominated by the region of finite solutions (gray area), indicating that no physically realizable folded configurations persist.
        This collapse can be understood by considering the integrable limit: as $N \to \infty$, the dynamics approaches a continuous integrable system where dissipation drives all trajectories toward the trivial attractor located at $I=0$.
        However, since the state $I=0$ corresponds to a geometric singularity (finite solution) in this origami model, trajectories inevitably hit the boundary of the admissible domain and terminate.
    
        In contrast, introducing a small drift term ($\mu = 0.05$, Fig.~\ref{fig:numerical_result_3}(B)) leads to a qualitative change in the phase-space structure.
        The drift term shifts the attractor away from the singular origin ($I=0$) to a location with a finite, non-zero action value.
        Consequently, a stable quasi-periodic attractor, which we refer to as a \emph{KAM attractor}, emerges within the physically valid domain.
        This result demonstrates that appropriate tuning of the drift parameter is essential for realizing robust origami structures in the presence of dissipation.
    
        Furthermore, by examining the evolution of the basins of attraction with increasing $N$, we identify a significant dynamical transition.
        As $N$ increases, the location of the point attractor gradually approaches the basin of the quasi-periodic attractor and eventually disappears.
        This behavior suggests that, in the continuum limit, the static equilibrium loses its distinct basin of attraction, leaving the quasi-periodic mode as the dominant stable state.
    
        Finally, we must address the geometric validity of the resulting structures. As in the symplectic case, the origami configurations corresponding to the KAM attractor are likely to exhibit self-intersections. Therefore, for physical fabrication, we recommend using fewer modules $N$.
        In this regime, although the system experiences stronger discrete perturbations, the lack of dense invariant curves facilitates the construction of tangible origami models.
    
\section{\label{sec:conclusion}CONCLUSION}
    In this study, we identified the inverse module number, $N^{-1}$, as a perturbation parameter and analyzed  tubular origami tessellations within the framework of KAM theory.

    Numerical experiments showed how the phase-space structure changes as this perturbation decreases.
    For sufficiently large module numbers, large regions of the phase space are occupied by structures that appear to be invariant curves, indicating that the finite-\(N\) system retains characteristic features of the integrable limit.
    Specifically, analysis of the frequency $\xi$ in the integrable limit proved effective for determining qualitative properties of the system, such as whether it exhibits twist or nontwist behavior, as well as for controlling the locations of stable folding regions corresponding to elliptic fixed points.
    
    A central theoretical contribution of this work is the geometric interpretation of the variational structure.
    We demonstrated that, in the integrable limit, the generating function of the folding map is identical to the total discrete mean curvature of the folded origami module.
    This result reveals a direct connection between the discrete geometry of folded surfaces and the variational structure underlying tubular origami dynamics.
    Moreover, the numerical results suggest that this curvature-based viewpoint remains useful for interpreting finite-\(N\) systems, even though the explicit form of the generating function becomes more complicated under perturbations.
    The resemblance between the folded configurations observed in this study and constant-mean-curvature surfaces such as \emph{Delaunay surfaces}\cite{Kenmotsu2003SurfacesCMC} further suggests a possible link between the origami dynamics and curvature-based geometry of folded surfaces.
    
    The phase-space description also provides a geometric interpretation of distinct stable folding modes.
    Previous work related quasiperiodic orbits around elliptic fixed points to undulated folded states in tubular origami tessellations.
    Building on this interpretation, continuous displacement between nearby orbits may be viewed as continuous deformation of the corresponding folded shapes.
    Our numerical computations indicate that invariant curves separating stable folding regions can be associated with self-intersecting folded shapes.
    Thus, these invariant curves can be interpreted as geometric barriers to continuous deformation between distinct stable folding modes.
    
    An important direction for future research is to extend the present analysis of the generating function, which was derived explicitly for symmetric crease configurations, to more general asymmetric geometries. In asymmetric systems, the variational structure is described using elliptic integrals, making its geometric interpretation analytically challenging. It remains an open question whether the correspondence between the generating function and curvature-based geometric quantities persists under perturbations, and how the integrable structure deforms in such generalized settings. 
    Furthermore, applying the insights from KAM theory obtained in this study to explore the design of infinite chaotic solutions is an important direction for future research. In the present system, there exist regions where solutions break down due to geometric constraints, resulting in finite solutions. Because typical chaotic orbits eventually enter these regions, the continuity of the folding dynamics is lost. However, by utilizing invariant curves as effective transport barriers in phase space, chaotic orbits can be confined within geometrically valid regions. This approach enables the design of infinitely continuous, non-periodic folded states while avoiding geometric failure. Exploiting the global structure of phase space to create origami tessellations corresponding to these chaotic solutions will provide deeper insights into the dynamical behavior of these systems.
    
    Finally, twist and nontwist behavior, as well as conformally symplectic structures arising from geometric degrees of freedom, may appear in broader classes of origami-based dynamical systems beyond the tubular setting.
    The integrable-limit approach developed in this study offers a useful perspective on the mathematical structures underlying complex origami dynamics. 
    Further examples and systematic comparisons will be needed to assess how broadly this viewpoint applies to origami engineering and metamaterial design.

\begin{acknowledgments}
M. S. was partially supported by JSPS KAKENHI Grant Number 23K25778.
We would like to thank Rinki Imada and Tomohiro Tachi, whose prior work on tubular origami tessellations informed this research, for providing the opportunity to discuss this study and for their valuable insights.
\end{acknowledgments}

\section*{Data Availability Statement}
Data sharing is not applicable to this article as no datasets were generated or analyzed.

\appendix

\section{\label{app:derivation}Appendix A: Analytical Derivation of the Generating Function}
    In this appendix, we provide a detailed algebraic derivation of the generating function
    $S_f(I)$ introduced in Sec.~\ref{sec:genfun}.
    The central task is the explicit evaluation of the integral term arising from the
    integration by parts.
    Substituting the geometric expressions for $A(I)$ and $B(I)$, the integrand
    $\psi(I)$ defined in Eq.~(\ref{eq:def_psi}) takes the form
    \begin{equation}
        \psi(I)
        =
        \frac{I(B'A - A'B)}{A^2 + B^2}
        =
        \frac{\sigma}{\beta\sqrt{\gamma}}\, R(I), \nonumber
    \end{equation}
    where $R(I)$ is a rational function of the action variable $I$.
    
    \subsection{\label{app:decomposition} Decomposition into Rational and Irrational Components}
        First, we express the integrand $\psi(I)$ explicitly in terms of the geometric parameters to separate its rational and irrational components.
        Observing the factor $I/(A^2+B^2)$, it is evident that this term is a rational function of $I$.
        Consequently, the algebraic nature of the integrand is determined by the term $B^{\prime}A - A^{\prime}B$.
        Recalling the definitions of $A(I)$ and $B(I)$ given in Eq.~(\ref{eq:AB_def}), these functions and their derivatives can be written as
        \begin{eqnarray}
            A &= \frac{R_{A}(I)}{\beta}, \qquad A^{\prime} = \frac{R_{A^{\prime}}(I)}{\beta}, \\
            B &= \sigma\sqrt{\gamma}, \qquad B^{\prime} = \frac{\sigma\gamma^{\prime}}{2\sqrt{\gamma}},
        \end{eqnarray}
        where $R_A(I)$ and $R_{A^{\prime}}(I)$ are rational functions given by
        \begin{eqnarray}
            R_A(I) &=& -e_1\beta^2 +\alpha\frac{r_1^2-r_3^2+4I^2-4Ie_1\alpha}{4I}, \\
            R_{A^{\prime}}(I) &=& R_{A}^{\prime}(I)+R_A(I)\frac{\alpha\alpha^{\prime}}{\beta^2}.
        \end{eqnarray}
        Substituting these expressions into the preceding term yields
        \begin{equation}
            B^{\prime}A - A^{\prime}B
            = \frac{\sigma}{\beta\sqrt{\gamma}}
            \left(\frac{\gamma^{\prime}}{2}R_A-\gamma R_{A^\prime}\right).
        \end{equation}
        Finally, by introducing the collective rational function
        \begin{equation}
            R(I) = \frac{I}{A^2+B^2}
            \left(\frac{\gamma^{\prime}}{2}R_A(I)-\gamma R_{A^\prime}(I)\right), \label{eq:def_R_app}
        \end{equation}
        we recover the factorized form given in Eq.~(\ref{eq:def_psi}), in which the irrational dependence is entirely captured by the factor $1/(\beta\sqrt{\gamma})$.
        
    \subsection{\label{app:polynomial_degree} Bifurcation of Integrability via the Degree of the Polynomial}
        Having isolated the rational component of the integrand, the analytic nature of the primitive function is entirely determined by the irrational denominator $\beta\sqrt{\gamma}$.
        By substituting the explicit geometric expressions for $\beta$ and $\gamma$ and simplifying the resulting algebraic structure, we find that this term reduces to the square root of a polynomial $P(I)$, normalized by a factor proportional to the action variable $I$:
        \begin{equation}
            \beta\sqrt{\gamma} = \sqrt{\beta^2\gamma} = \frac{\sqrt{P(I)}}{8 I l_L^2}.
        \end{equation}
        Here, $P(I)=p_0+p_1+p_2$ is a polynomial in $I$, whose coefficients are fully determined by the geometric parameters of the module.
        The individual terms $p_j \ (j=0,1,2)$ are given explicitly as follows:
        \begin{eqnarray}
            p_0 &=& 4l_{L}^{2} r_{1}^{2} \left(16 I^{2} l_{L}^{2} - \left(4 I^{2} + l_{L}^{2} - l_{R}^{2}\right)^{2}\right), \\
            p_1 &=& \left(-16 I^{2} l_{L}^{2} + \left(4 I^{2} + l_{L}^{2} - l_{R}^{2}\right)^{2}\right) \left(l_{L}^{2} + r_{1}^{2} - r_{2}^{2}\right)^{2}, \\
            p_2 &=& -\bigg( 2 l_{L}^{2} \left(4 I^{2} + r_{1}^{2} - r_{3}^{2}\right) \nonumber \\
                & & \quad - \left(4 I^{2} + l_{L}^{2} - l_{R}^{2}\right) \left(l_{L}^{2} + r_{1}^{2} - r_{2}^{2}\right) \bigg)^{2}.
        \end{eqnarray}
        From the explicit form of $P(I)$, it follows that in the general asymmetric configuration ($l_L \neq l_R$ and $r_1 \neq r_3$), the highest-order terms in $I$ do not cancel.
        Consequently, the polynomial $P(I)$ is of degree four (quartic), implying that the primitive function associated with $\psi(I)$ belongs to the class of elliptic integrals.
        
        In contrast, for the symmetric configuration ($l_L = l_R = l$ and $r_1 = r_3$), the algebraic structure simplifies drastically.
        In this case, the polynomial $P(I)$ factorizes, allowing the cancellation of the term $I$ in the denominator.
        Consequently, the irrational factor reduces to the following form:
        \begin{equation}
            \beta\sqrt{\gamma}
            =
            \frac{\sqrt{K_1 - K_2 I^2}}{2l},
        \end{equation}
        where the constants $K_1$ and $K_2$ depend solely on the geometric parameters and are given by:
        \begin{eqnarray}
            K_1 &=& -l^{4} + 2 l^{2} (r_1^{2} + r_{2}^{2}) - (r_1^{2} - r_{2}^{2})^2, \\
            K_2 &=& 4r_2^2.
        \end{eqnarray}
        As a result, the term inside the square root becomes a quadratic function of $I$.
        This reduction in degree ensures that the integral of $\psi(I)$ can be evaluated explicitly in terms of elementary functions.

    \subsection{\label{app:integration_result} Explicit Form of the Integration Result}
        As shown in the previous subsection, imposing the symmetry conditions
        ($l_L = l_R = l$ and $r_1 = r_3$) reduces the integrand $\psi(I)$ to a form whose
        primitive can be expressed in terms of elementary functions.
        We briefly summarize the evaluation of this integral.
        
        Introducing the substitution
        \begin{equation}
            I = K \sin \rho, \qquad K = \sqrt{\frac{K_1}{K_2}},
        \end{equation}
        followed by the \emph{Weierstrass substitution}
        \begin{equation}
            t = \tan\left(\frac{\rho}{2}\right),
        \end{equation}
        the integral simplifies to a rational function of $t$.
        Specifically, one obtains
        \begin{equation}
            \int \psi(I)\, dI
            =
            \frac{\sigma l}{r_2}
            \int
            R\!\left(\frac{2K t}{1+t^2}\right)
            \frac{2}{1+t^2}
            \, dt,
        \end{equation}
        where $R(I)$ is the rational function defined in
        Eq.~(\ref{eq:def_R_app}).
        The resulting integral can be evaluated by standard partial fraction
        decomposition. Conducting the decomposition explicitly, the integrand
        can be organized into three characteristic contributions associated with
        the geometric parameters $l$, $r_1$, and $r_2$, respectively:

        \begin{widetext}
            \begin{eqnarray}
                & &\int \psi(I)\, dI = \nonumber \\
                &=&  \frac{\sigma l}{r_2}\int 
                        \frac{
                            4 K^{2}  t^{2}
                            \left(
                                8 K^{2} r_{2}^{2} t^{2}
                                + (t^{2} + 1)^{2}
                                  \left(
                                      l^{4}
                                      - 2 l^{2} r_{1}^{2}
                                      - l^{2} r_{2}^{2}
                                      + r_{1}^{4}
                                      - r_{1}^{2} r_{2}^{2}
                                  \right)
                            \right)
                        }{
                            l (t^{2} + 1)
                            \left(
                                16 K^{4} t^{4}
                                - 4 K^{2} t^{2} (l^{2} + r_{1}^{2}) (t^{2} + 1)^{2}
                                + l^{2} r_{1}^{2} (t^{2} + 1)^{4}
                            \right)
                        } \ dt, \nonumber \\
                &=& \frac{\sigma l}{r_2}\int \Biggl(
                    \frac{\left(-l^{2} + r_{1}^{2} - r_{2}^{2}\right)}
                         {2 \left(2 K t + l t^{2} + l\right)}
                    + \frac{\left(-l^{2} + r_{1}^{2} - r_{2}^{2}\right)}
                           {2 \left(-2 K t + l t^{2} + l\right)}
                    - \frac{r_{1} \left(-l^{2} + r_{1}^{2} + r_{2}^{2}\right)}
                           {2 l \left(2 K t + r_{1} t^{2} + r_{1}\right)}
                    - \frac{r_{1} \left(-l^{2} + r_{1}^{2} + r_{2}^{2}\right)}
                           {2 l \left(-2 K t + r_{1} t^{2} + r_{1}\right)}
                    + \frac{2 r_{2}^{2}}
                           {l \left(t^{2} + 1\right)}
                    \Biggr) \ dt, \nonumber \\
                &=& 
                \frac{\sigma l}{r_2}\Bigg(
                 \frac{-l^2+r_1^2-r_2^2}{2\sqrt{l^2-K^2}}
                     \left(
                        \tan^{-1}\left(
                        \frac{lt+K}{\sqrt{l^2-K^2}}
                        \right) + 
                         \tan^{-1}\left(
                        \frac{lt-K}{\sqrt{l^2-K^2}}
                        \right)
                     \right) \nonumber \\ 
                 & & \qquad \qquad -
                     \frac{r_1}{l}
                     \frac{-l^2+r_1^2+r_2^2}{2\sqrt{r_1^2-K^2}}
                     \left(
                        \tan^{-1}\left(
                        \frac{r_1t+K}{\sqrt{r_1^2-K^2}}
                        \right) + 
                         \tan^{-1}\left(
                        \frac{r_1t-K}{\sqrt{r_1^2-K^2}}
                        \right)
                     \right)
                     + \frac{2r_2^2}{l}\tan^{-1}t 
                 \Bigg) + \mathrm{Const} . \label{eq:long_integral}
            \end{eqnarray}
        \end{widetext}
        Although the integrated expression in Eq.~(\ref{eq:long_integral}) appears algebraically complex, it can be significantly simplified by considering the geometric interpretation of each term.
        \begin{figure}
            \includegraphics[width=\columnwidth]{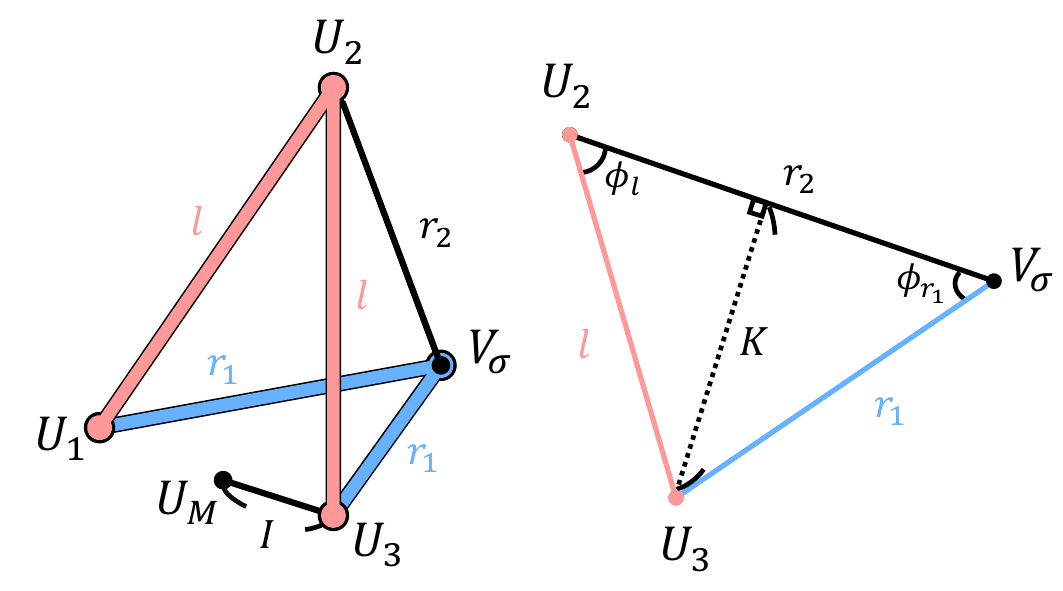}
            \caption{\label{fig:app_triangle_with_K}
                 Geometric interpretation of the auxiliary parameter $K$ and the partial angles on the polyhedral faces. The constant $K$ represents the altitude of the triangle with base $r_2$, while $\phi_l$ and $\phi_{r_1}$ denote the angles obtained from the transformation of the inverse tangent terms.
            }
        \end{figure}
        By invoking the \emph{Law of Cosines} and \emph{Heron's formula}, the coefficient $K$ introduced in the variable transformation is identified as the altitude of the triangle $U_3 U_2 V_\sigma$ with respect to the base $r_2$ (see Fig.~\ref{fig:app_triangle_with_K}).
        Furthermore, by defining the auxiliary angles $\phi_l$ and $\phi_{r_1}$ as depicted in the figure, the integrated expression allows for a significant simplification.
        Specifically, the coefficients reduce to:
        \begin{eqnarray}
            \frac{-l^2+r_1^2-r_2^2}{2\sqrt{l^2-K^2}} &=& -r_2, \\
            -\frac{r_1}{l}\frac{-l^2+r_1^2+r_2^2}{2\sqrt{r_1^2-K^2}} &=& -\frac{r_1r_2}{l} ,
        \end{eqnarray}
        and the inverse trigonometric terms are rewritten using the addition formulas as:
        \begin{eqnarray}
            \Theta_l &:=& \tan^{-1}\left(
            \frac{lt+K}{\sqrt{l^2-K^2}}
            \right) + 
             \tan^{-1}\left(
            \frac{lt-K}{\sqrt{l^2-K^2}}
            \right) \nonumber\\
            &=& \tan^{-1}\left(\frac{2t}{1-t^2}\cos\phi_l\right), \\
            \Theta_{r_1} &:=&  \tan^{-1}\left(
                \frac{r_1t+K}{\sqrt{r_1^2-K^2}}
                \right) + 
                 \tan^{-1}\left(
                \frac{r_1t-K}{\sqrt{r_1^2-K^2}}
                \right) \nonumber\\
            &=& \tan^{-1}\left(\frac{2t}{1-t^2}\cos\phi_{r_1}\right).
        \end{eqnarray}
        Recalling the definition of the Weierstrass substitution $t = \tan(\phi/2)$, we identify the rational term involving $t$ as the tangent of the angle $\rho$:
        \begin{equation}
            \frac{2t}{1-t^2} = \tan\rho.
        \end{equation}
        Consequently, all algebraic complexity reduces to a notably simple geometric form:
        \begin{equation}
            -\int \psi(I) \, dI = \sigma\left(l\Theta_l + r_1\Theta_{r_1} - r_2\rho\right) + \mathrm{Const.}
        \end{equation}
        
        \subsection{\label{app:geometric_interpretation} Geometric Interpretation of the Generating Function}
            Combining the explicit integration result of $\psi(I)$ derived in the previous section, we now give a geometric interpretation of the generating function $S_f$ (Eq.~(\ref{eq:S_f_first_def})).
            We first clarify the meaning of $\xi(I)$, which is defined through the auxiliary functions $A(I)$ and $B(I)$ (see Eqs.~(\ref{eq:AB_def})).
            
\begin{figure}
                \includegraphics[width=\columnwidth]{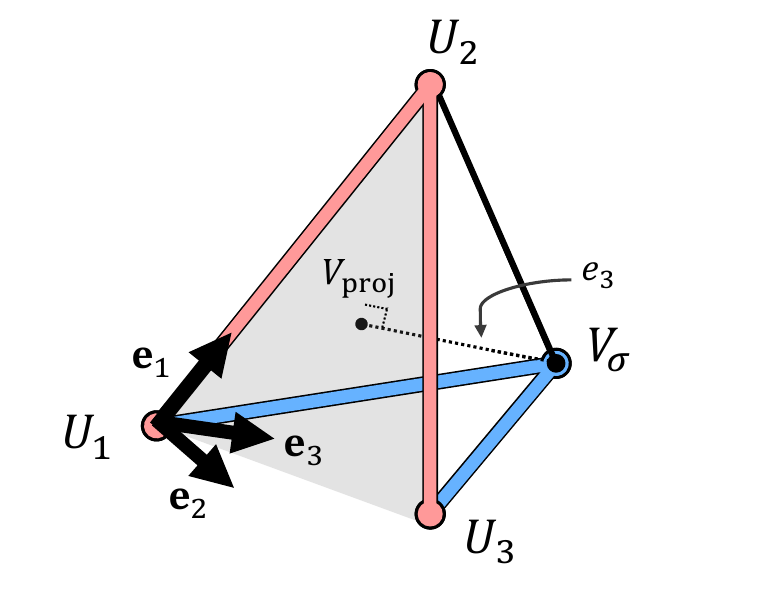}
                \caption{\label{fig:app_AB_meaning}
                    Projection geometry defining $A(I)$ and $B(I)$.
                    $V_{\mathrm{proj}}$ denotes the orthogonal projection of $V_\sigma$ onto the plane
                    $U_1U_2U_3$.
                }
            \end{figure}

            Following the model setup in Sec.~\ref{sec:model}, we introduce a local coordinate system
            $(\mathbf{e}_1,\mathbf{e}_2,\mathbf{e}_3)$ centered at the vertex $U_1$
            (Fig.~\ref{fig:app_AB_meaning}).
            In this frame, $B(I)$ is the $\mathbf{e}_3$ component and represents the signed perpendicular distance from $V_\sigma$ to the plane $U_1U_2U_3$.
            Let $V_{\mathrm{proj}}$ have coordinates $[e_1,e_2]^T$ in the $(\mathbf{e}_1,\mathbf{e}_2)$ basis.
            With $\alpha=\cos\angle U_2U_1U_3$ and $\beta=\sin\angle U_2U_1U_3$, we obtain
            \begin{equation}
                \begin{bmatrix}
                    \alpha & \beta \\
                    -\beta & \alpha
                \end{bmatrix}
                \begin{bmatrix}
                    e_1 \\ e_2
                \end{bmatrix}
                =
                \begin{bmatrix}
                    \cdots \\ A(I)
                \end{bmatrix}.
            \end{equation}
            This rotation aligns $\mathbf{e}_1$ with the edge $U_1U_3$.
            Hence, $A(I)$ corresponds (up to sign) to the transverse component in the rotated frame.
            Together, $A(I)$ and $B(I)$ define the exterior dihedral angle $\xi(I)$ associated with the edge $U_1U_3$.
            
\begin{figure}
                \includegraphics[width=\columnwidth]{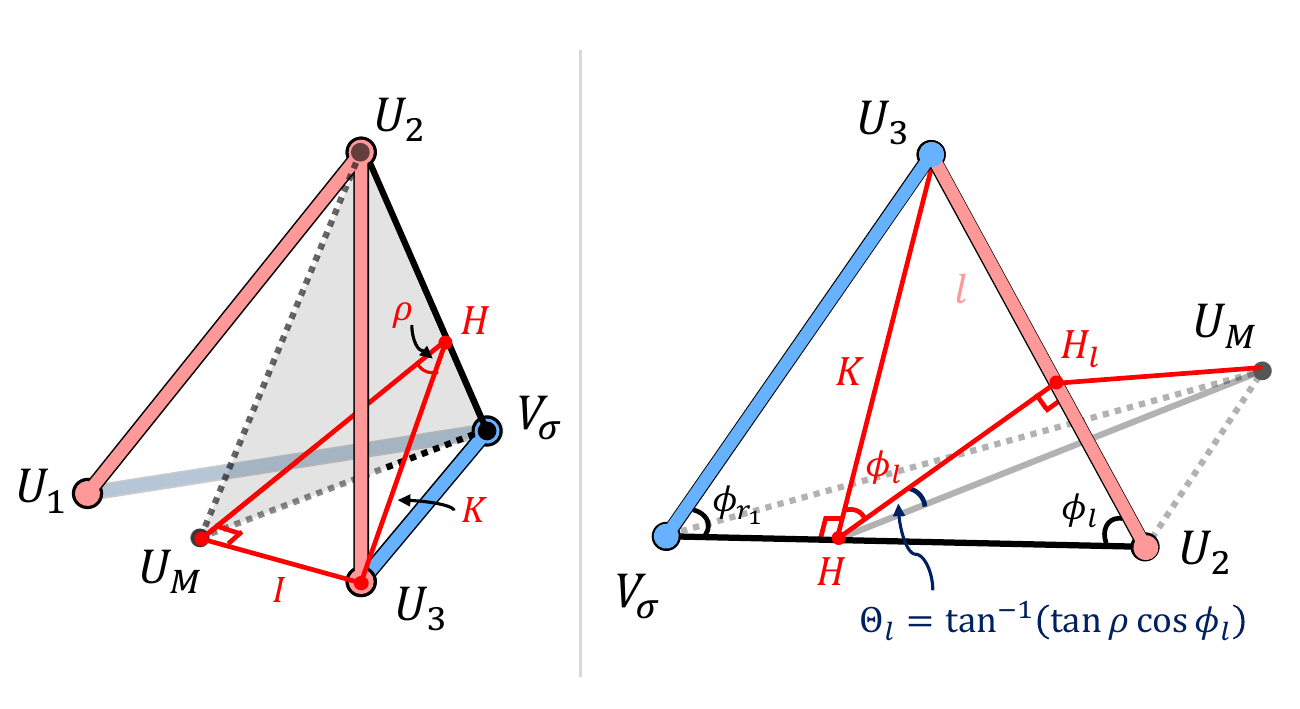}
                \caption{\label{fig:app_angles}
                    (Left) Definition of $\rho$ via $I=K\sin\rho$.
                    (Right) Construction of the angle $\Theta_l$ associated with the edge $l$.
                }
            \end{figure}

            We now interpret the angular variables $\Theta_l$, $\Theta_{r_1}$, and $\rho$ associated with the edges
            $l$, $r_1$, and $r_2$, respectively.
            The substitution $I=K\sin\rho$ identifies $\rho$ as half of the interior dihedral angle along $r_2$
            (Fig.~\ref{fig:app_angles}, left).
            The integrated term $-r_2\rho$ is therefore equivalent, up to an additive constant,
            to the contribution of the exterior dihedral angle associated with the edge $r_2$,
            which we denote by $\Phi_{r_2}$.
            It can be written as $r_2\Phi_{r_2}/2$, and such constant shifts are irrelevant
            in variational calculations.
            
            For the edge $l$, let $H_l$ be the foot of the perpendicular from $H$ onto $U_2U_3$
            (Fig.~\ref{fig:app_angles}, right), and let $\Pi_1$ and $\Pi_2$ be the planes $U_2U_3V_\sigma$ and $U_2U_MV_\sigma$, respectively.
            The dihedral angle between $\Pi_1$ and $\Pi_2$ is $\rho$.
            If the direction $H_lH$ is rotated by $\phi_l$ from the direction of maximal slope,
            its effective inclination satisfies
            \begin{equation}
                \tan\Theta_l=\tan\rho\,\cos\phi_l .
            \end{equation}
            
            From the planar geometry of the triangle $H_lHU_M$, the exterior dihedral angle $\Phi_l$ obeys
            \begin{equation}
                \Theta_l + \frac{\pi}{2} + (\pi - \Phi_l) = \pi.
            \end{equation}
            Since the difference is a constant, $\Theta_l$ is variationally equivalent to $\Phi_l$.
            An analogous construction shows that $\Theta_{r_1}$ is variationally equivalent to
            the exterior dihedral angle $\Phi_{r_1}$.
                      
\section{\label{app:schlafli_derivation}Derivation of the Mapping Equations from the Explicit Generating Function}
    In this appendix, we verify that the explicit geometric form of the generating function $S_f$ derived above correctly reproduces the original mapping equations.
    Specifically, this verification is performed by applying the canonical relations established in Eq.~(\ref{eq:canonical_relationship}) and the frequency definition in Eq.~(\ref{eq:gf_corresponds_xi}).
    
    Recall that $S_f$ represents the total discrete mean curvature of the tetrahedron $U_1 U_2 U_3 \text{-} V_\sigma$, expressed as the sum of the products of edge lengths and their corresponding exterior dihedral angles.
    For simplicity, let us denote the lengths of the six edges of the tetrahedron by
    $L_k$ and their corresponding exterior dihedral angles by $\Phi_k$
    ($k = 1, \dots, 6$).
    In particular, we assign the index $k=1$ to the diagonal edge $U_1U_3$,
    whose length is given by $L_1 = 2I$,
    while the remaining indices correspond to the fixed edges of the crease pattern.
    The generating function $S_f$ is thus written as:
    \begin{equation}
        S_f = \frac{1}{2}\sum_{k=1}^6 L_k \Phi_k.
    \end{equation}
    Differentiating this expression with respect to the action variable $I$ yields two distinct terms:
    \begin{equation}
        \frac{d S_f}{d I} = \frac{1}{2}\sum_{k=1}^6 \Phi_k \frac{d L_k}{d I} + \frac{1}{2}\sum_{k=1}^6 L_k \frac{d \Phi_k}{d I}.
    \end{equation}
    The second term on the right-hand side represents the variation of the dihedral angles, weighted by the edge lengths.
    According to the \emph{Schläfli formula}\cite{milnor1994schlafli, regge1961general} for a closed polyhedron (specifically, a tetrahedron), this sum is identically zero:
    \begin{equation}
        \sum_{k=1}^6 L_k \, d\Phi_k = 0 \quad \Longrightarrow \quad \sum_{k=1}^6 L_k \frac{d \Phi_k}{dI} = 0.
    \end{equation}
    Consequently, only the first term remains.
    In the present model, the edge lengths associated with the fixed crease pattern,
    namely $l$, $r_1$, and $r_2$, are constant.
    The only edge whose length depends on $I$ is the diagonal edge $U_1U_3$,
    whose length is given by $L_{1} = 2I$.
    Let $\Phi_{1} = \xi_f(I)$ denote the exterior dihedral angle associated with this edge. The derivative of the generating function therefore simplifies to
    \begin{equation}
        \frac{d S_f}{dI} = \xi_f(I).
    \end{equation}
    Thus, differentiating the geometric generating function $S_f$ with respect to the action variable yields the corresponding conjugate angular variable.

\nocite{*}
\bibliography{aipsamp}

\end{document}